\documentclass[11pt, reqno]{amsart}

\usepackage[usenames,dvipsnames]{color}

\usepackage[OT2, T1]{fontenc}
\usepackage{url}
\usepackage{amsmath}
\usepackage{array}
\usepackage{graphicx}
\usepackage{amsfonts}
\usepackage{amssymb}
\usepackage{amsthm}
\usepackage{hyperref}
\usepackage{paralist}
\usepackage{colonequals}	
\usepackage{color}
\usepackage{mathabx}		
\usepackage{amscd}
\usepackage[all,cmtip]{xy}	
\usepackage{sseq}			
\usepackage{verbatim}
\usepackage{parskip}
\usepackage{microtype}
\usepackage[in]{fullpage}
\usepackage{mathrsfs} 			
\usepackage{bm} 				
\usepackage{cleveref}
\usepackage{enumitem}
\usepackage{moreenum}			
\usepackage[alphabetic,lite]{amsrefs} 	
\usepackage{stmaryrd}	

\DeclareSymbolFont{cyrletters}{OT2}{wncyr}{m}{n}
\DeclareMathSymbol{\Sha}{\mathalpha}{cyrletters}{"58}

\newcommand{\gA}{\alpha}

\newcommand{\bA}{\mathbb{A}}

\newcommand{\bF}{\mathbb{F}}
\newcommand{\bG}{\mathbb{G}}

\newcommand{\bQ}{\mathbb{Q}}
\newcommand{\bR}{\mathbb{R}}

\newcommand{\bZ}{\mathbb{Z}}

\newcommand{\bbB}{\mathbf{B}}

\newcommand{\bbR}{\mathbf{R}}

\newcommand{\cC}{\mathcal{C}}

\newcommand{\cF}{\mathcal{F}}
\newcommand{\cG}{\mathcal{G}}
\newcommand{\cH}{\mathcal{H}}

\newcommand{\cM}{\mathcal{M}}

\newcommand{\cO}{\mathcal{O}}

\newcommand{\cQ}{\mathcal{Q}}

\newcommand{\cT}{\mathcal{T}}

\newcommand{\fo}{\mathfrak{o}}

\newcommand{\ra}{\rightarrow}

\newcommand{\xra}{\xrightarrow}

\newcommand{\hra}{\hookrightarrow}

\newcommand{\wt}{\widetilde}
\newcommand{\wh}{\widehat}

\newcommand{\pr}{^{\prime}}

\newcommand{\ce}{\colonequals}
\newcommand{\ov}{\overline}

\renewcommand{\b}{\textbf}
\newcommand{\surjects}{\twoheadrightarrow}

\newcommand{\tensor}{\otimes} 		
\newcommand{\isomto}{\overset{\sim}{\longrightarrow}}

\newcommand{\nr}{{\mathrm{nr}}}		
\newcommand{\cont}{{\mathrm{cont}}}		
\newcommand{\et}{\mathrm{\acute{e}t}}	
\newcommand{\tors}{\mathrm{tors}}		
\newcommand{\sh}{\mathrm{sh}}		

\providecommand{\abs}[1]{\left\lvert#1\right\rvert}

\providecommand{\p}[1]{\left(#1\right)}

\providecommand{\f}[2]{\frac{#1}{#2}}
\DeclareMathOperator{\Ker}{Ker}			
\DeclareMathOperator{\im}{Im}			
\DeclareMathOperator{\Spec}{Spec}		
\DeclareMathOperator{\Hom}{Hom}			
\DeclareMathOperator{\Char}{char}		
\DeclareMathOperator{\loc}{loc}		%
\DeclareMathOperator{\Br}{Br}		
\DeclareMathOperator{\Gal}{Gal}	
\DeclareMathOperator{\inv}{inv}	
\DeclareMathOperator{\Res}{Res}		
\DeclareMathOperator{\GL}{GL}		
\DeclareMathOperator{\Aut}{Aut}		
\DeclareMathOperator{\Sel}{Sel}		
\DeclareMathOperator{\Pic}{Pic}		
\newcommand{\ba}{\begin{aligned}}
\newcommand{\ea}{\end{aligned}}
\newcommand{\be}{\begin{equation}}
\newcommand{\ee}{\end{equation}}
\newcommand{\pf}{\begin{proof}}
\newcommand{\bpf}{\begin{proof}}
\newcommand{\epf}{\end{proof}}
\newcommand{\bthm}{\begin{thm}}
\newcommand{\ethm}{\end{thm}}
\newcommand{\bprop}{\begin{prop}}
\newcommand{\eprop}{\end{prop}}
\newcommand{\bcor}{\begin{cor}}
\newcommand{\ecor}{\end{cor}}
\newcommand{\brem}{\begin{rem}}
\newcommand{\erem}{\end{rem}}
\newcommand{\brems}{\begin{rems} \hfill \begin{enumerate}[label=\b{\thesubsection.},ref=\thesubsection]}
\newcommand{\remi}{\addtocounter{subsection}{1} \item}
\newcommand{\erems}{\end{enumerate} \end{rems}}
\newcommand{\blem}{\begin{lemma}}
\newcommand{\elem}{\end{lemma}}
\newcommand{\bconj}{\begin{conj}}
\newcommand{\econj}{\end{conj}}
\newcommand{\bprob}{\begin{Problem}}
\newcommand{\eprob}{\end{Problem}}
\newcommand{\bq}{\begin{q}}
\newcommand{\eq}{\end{q}}
\newcommand{\benum}{\begin{enumerate}[label={(\alph*)}]}
\newcommand{\benuma}{\begin{enumerate}[label={(\arabic*)}]}
\newcommand{\benumr}{\begin{enumerate}[label={(\roman*)}]}
\newcommand{\eenum}{\end{enumerate}}
\newcommand{\bc}{}
\newcommand{\beg}{\begin{eg}}
\newcommand{\eeg}{\end{eg}}
\newcommand{\bcl}{\begin{claim}}
\newcommand{\ecl}{\end{claim}}
\newcommand{\lab}{\label}
\newcommand{\tst}{\textstyle}
\newcommand{\qq}{\quad\quad}
\newcommand{\qqq}{\quad\quad\quad}

\newcommand*{\QED}{\hfill\ensuremath{\qed}}
\newcommand*{\QEDD}{\hfill\ensuremath{\qed\qed}}

\theoremstyle{plain}
\newtheorem{thm}[subsection]{Theorem}
\Crefname{thm}{Theorem}{Theorems}

\Crefname{rethm}{Theorem}{Theorem}
\newtheorem{prop}[subsection]{Proposition}
\Crefname{prop}{Proposition}{Propositions}
\newtheorem{q}[subsection]{Question}
\Crefname{q}{Question}{Questions}
\Crefname{eg}{Example}{Examples}
\newtheorem{Problem}[subsection]{Problem}
\Crefname{Problem}{Problem}{Problems}
\newtheorem{conj}[subsection]{Conjecture}
\Crefname{conj}{Conjecture}{Conjectures}
\newtheorem{cor}[subsection]{Corollary}
\Crefname{cor}{Corollary}{Corollaries}

\newtheorem{lemma}[subsection]{Lemma}

\Crefname{subprop}{Proposition}{Propositions}
\Crefname{subcor}{Corollary}{Corollaries}
\Crefname{sublem}{Lemma}{Lemmas}

\theoremstyle{remark}
\newtheorem{claim}[equation]{Claim}
\Crefname{claim}{Claim}{Claims}

\Crefname{subrem}{Remark}{Remarks}

\theoremstyle{definition}
\newtheorem{eg}[subsection]{Example}
\newtheorem{rem}[subsection]{Remark}
\Crefname{rem}{Remark}{Remarks}
\newtheorem*{rems}{Remarks}

\newtheoremstyle{subsection-tweak}
   {11pt}
   {3pt}%
   {}
   {}%
   {\bfseries}
   {}%
   {.5em}
   {\thmnumber{\@{#1}{}\@{#2}.}%
    \thmnote{~{\bfseries#3.}}}

\Crefname{innercustomconj}{Conjecture}{Conjecture}

\theoremstyle{subsection-tweak}
\newtheorem{pp}[subsection]{}
\newcommand{\bpp}{\begin{pp}}
\newcommand{\epp}{\end{pp}}

\numberwithin{equation}{subsection}

\begin{document}
\author{K\k{e}stutis \v{C}esnavi\v{c}ius}
\title{Poitou--Tate without restrictions on the order}
\date{\today}
\address{Department of Mathematics, University of California, Berkeley, CA 94720-3840, USA}
\email{kestutis@berkeley.edu}
\urladdr{http://math.berkeley.edu/~kestutis/}

\begin{abstract} 
The Poitou--Tate sequence relates Galois cohomology with restricted ramification of a finite Galois module $M$ over a global field to that of the dual module under the assumption that $\#M$ is a unit away from the allowed ramification set. We remove the assumption on $\#M$ by proving a generalization that allows arbitrary ``ramification sets'' that contain the archimedean places. We also prove that restricted products of local cohomologies that appear in the Poitou--Tate sequence may be identified with derived functor cohomology of an adele ring. In our proof of the generalized sequence we adopt this derived functor point of view and exploit properties of a natural topology carried by cohomology of the adeles. 
  \end{abstract}

\maketitle

\section{Introduction}

\bpp[Poitou--Tate for Galois cohomology with restricted ramification] \lab{PT-Gal}
One of the core results of global arithmetic duality is the Poitou--Tate exact sequence of Galois cohomology groups:
\be\tag{$\star$}\lab{PT-old}\ba
\xymatrix@C=20pt @R=12pt{
  0  \ar[r] &   H^0(\Gal(K_{\Sigma}/K), M) \ar[r] & \prod_{v\in \Sigma}\pr H^0(K_v, M)   \ar[r] & H^2(\Gal(K_{\Sigma}/K), M^D)^*               \ar@{->} `r/8pt[d] `/9pt[l] `^dl/9pt[lll] `^r/7pt[dl] [dll] \\
&  H^1(\Gal(K_{\Sigma}/K), M) \ar[r] & \prod_{v\in \Sigma}\pr H^1(K_v, M)   \ar[r] & H^1(\Gal(K_{\Sigma}/K), M^D)^*        \ar@{->} `r/8pt[d] `/9pt[l] `^dl/9pt[lll] `^r/7pt[dl] [dll] \\
&  H^2(\Gal(K_{\Sigma}/K), M) \ar[r] & \prod_{v\in \Sigma}\pr H^2(K_v, M)   \ar[r] & H^0(\Gal(K_{\Sigma}/K), M^D)^*  \ar[r] & 0,
}
\ea\ee
where $K$ is a global field, $\Sigma$ is a possibly infinite nonempty set of places of $K$ containing the archimedean places, $\Gal(K_{\Sigma}/K)$ is the Galois group of the maximal unramified outside $\Sigma$ subextension $K_{\Sigma}/K$ of a separable closure $K^s/K$ (a choice of which is fixed), $M$ is a finite discrete $\Gal(K_\Sigma/K)$-module whose order is a unit at all places not in $\Sigma$, the restricted products $\prod_{v\in \Sigma}\pr$ are with respect to the unramified subgroups, and $M^D\ce \Hom_\bZ(M, \bQ/\bZ(1))$ is the dual of $M$.

The Poitou--Tate sequence \eqref{PT-old} is proved in \cite{Mil06}*{I.4.10} and is useful in many contexts. However, its requirement that $\#M$ be a unit outside $\Sigma$ may be restrictive. The goal of this paper is to dispose of this requirement by proving a generalization of \eqref{PT-old} and to offer another, possibly more conceptual, interpretation of the restricted products appearing in \eqref{PT-old}. Our generalization allows new $\Sigma$ even in the number field case, and it also gives a different proof of the exactness of \eqref{PT-old}.
\epp

\bpp[Poitou--Tate without restrictions on $\#M$] \lab{PT-wo}
As in \S\ref{PT-Gal}, let $\Sigma$ be a set of places containing the archimedean places. For simplicity, suppose that $\Sigma$ also contains a finite place (see \Cref{PT-thm} for the general case). The sought generalization of \eqref{PT-old} uses fppf cohomology and reads
\be\tag{$\bigstar$}\lab{PT-new}\ba
\xymatrix@C=20pt @R=12pt{
   0  \ar[r] &   H^0(S - \Sigma, G) \ar[r] & \wh{H}^0(\bA_K^{\in \Sigma}, G)   \ar[r] & H^2(S - \Sigma, H)^*               \ar@{->} `r/8pt[d] `/9pt[l] `^dl/9pt[lll] `^r/7pt[dl] [dll] \\
&  H^1(S - \Sigma, G) \ar[r] & H^1(\bA_K^{\in \Sigma}, G)   \ar[r] & H^1(S - \Sigma, H)^*        \ar@{->} `r/8pt[d] `/9pt[l] `^dl/9pt[lll] `^r/7pt[dl] [dll] \\
&  H^2(S - \Sigma, G) \ar[r] & H^2(\bA_K^{\in \Sigma}, G)   \ar[r] & H^0(S - \Sigma, H)^*  \ar[r] & 0, 
}
\ea\ee
where, as defined in \S\ref{not}, $S - \Sigma$ is the localization away from $\Sigma$ of the smooth proper curve with function field $K$ if $\Char K > 0$ (resp.,~of the spectrum of the ring of integers of $K$ if $\Char K = 0$), $G$ and $H$ are Cartier dual commutative finite flat $(S - \Sigma)$-group schemes, $\bA_K^{\in\Sigma}$ is the adele ring that takes only the places in $\Sigma$ into account, $H^n(\bA_K^{\in\Sigma}, -)$ is its fppf cohomology, and $\wh{H}^0(\bA_K^{\in\Sigma}, -)$ is $H^0(\bA_K^{\in\Sigma}, -)$ unless $K$ is a number field (in which case, see \S\ref{H0-Tate} for the definition).\footnote{In \eqref{PT-old}, the terms $H^0(K_v, M)$ for archimedean $v$ are also not the honest $H^0$ but their Tate modifications.}

As we explain in \Cref{app}, if $M$ of \S\ref{PT-Gal} has order that is a unit outside $\Sigma$, then $M$ and $M^D$ extend to Cartier dual commutative finite \'{e}tale $(S - \Sigma)$-group schemes whose fppf cohomology groups identify with Galois cohomology groups used in \eqref{PT-old} and the sequence \eqref{PT-new} recovers \eqref{PT-old}.
\epp

\bpp[Restricted products in other contexts]
Although the adelic interpretation of the restricted products is an aesthetic improvement to \eqref{PT-old}, it is an avatar of a general phenomenon observed in \Cref{res-pr-comp}: if $G$ is a $K$-group scheme of finite type, then $H^1(\bA_K, G)$ decomposes as a restricted product of local cohomologies whenever $G$ is either smooth connected or proper, and so does $H^n(\bA_K, G)$ with $n \ge 2$ whenever $G$ is commutative. If one wishes, one may use this to reinterpret restricted products appearing elsewhere, e.g.,~as explained in \Cref{CFT}, one may recast the reciprocity sequence of global class field theory as an exact sequence relating the Brauer group of $K$ to that of $\bA_K$.
\epp

\brems 
\remi \lab{GA-rem}
In the special case when $\Char K > 0$ and $\Sigma$ contains all places, the sequence \eqref{PT-new} (without the adelic interpretation of the restricted products) was proved in \cite{GA09}*{4.11}. We follow the strategy used there but need additional input---the proof of loc.~cit.~could only give \eqref{PT-new} after enlarging $\Sigma$ by finitely many places that depend on $G$. The new input is provided by \cite{Ces14a}*{\S2}, which improves several auxiliary results of \cite{GA09} by ensuring that no enlargement of $\Sigma$ is needed.

\remi
The proof of \eqref{PT-old} given in \cite{Mil06}*{I.4.10} seems to be close in spirit to that of global class field theory, being rooted in manipulations of Galois cohomology of various modules of arithmetic interest. The proof of \eqref{PT-old} that results as a special case of \eqref{PT-new} is based on geometric ideas: its main global input is the generalization \cite{Mil06}*{III.3.1 and III.8.2} of Artin--Verdier duality. The method of deriving other arithmetic duality theorems from the Artin--Verdier duality has also been used elsewhere, see, for instance, \cite{Mil06}*{II.\S5}, \cite{HS05}, or \cite{GA09}.
\erems

\bpp[An overview of the proof and of the paper]
We begin in \S\ref{coho-res} by comparing derived functor cohomology of adele rings with  restricted products of local cohomologies. We use this comparison in \S\ref{topo-ad} to show that cohomology of the adeles carries a canonical topology; we then proceed to record basic properties of this topology. Here the results of \cite{Ces15d} come in handy by supplying the corresponding topological properties in a local setting. 

Topological aspects, and hence also our work in \S\ref{topo-ad}, are crucial for the proof of the exactness of \eqref{PT-new}. However, a reader who is willing to accept the necessary topological properties of $H^n(\bA_K^{\in \Sigma}, -)$ as known may skip \S\ref{topo-ad} and instead consult \S\ref{basic-topo}, which gathers together those conclusions of \S\S\ref{coho-res}--\ref{topo-ad} that are needed in \S\S\ref{disc-image}--\ref{CPT}. Likewise, a reader who is interested in the proof of \eqref{PT-new} but would prefer to take $H^n(\bA_K^{\in\Sigma}, -)$ to mean $\prod_{v\in\Sigma}\pr H^n(K_v, -)$ may skip \S\S\ref{coho-res}--\ref{topo-ad} entirely and start with \S\ref{disc-image}.

In \S\ref{disc-image} we prove arithmetic duality results that are key for the proof of \eqref{PT-new}; \S\ref{disc-image} is also the place where the input from \cite{Ces14a}*{\S2} mentioned in Remark \ref{GA-rem} enters the picture (through \Cref{Sha}). In \S\ref{PT} we show that the sequence \eqref{PT-new} is a quick consequence of the results of \S\ref{disc-image}. In \S\ref{CPT} we include a proof of the Cassels--Poitou--Tate sequence in the generality permitted by \eqref{PT-new}. This sequence with finite $\Sigma$ has played instrumental roles in the proofs of a number of results on the arithmetic of abelian varieties, so it seems worthwhile to record in \S\ref{CPT} its generalization that allows arbitrary $\Sigma$. We conclude the paper by showing in \Cref{app} that \eqref{PT-new} indeed generalizes \eqref{PT-old}.
\epp

\bpp[Notation] \lab{not}
The following notation is in place throughout the paper: $K$ is a global field; if $\Char K = 0$, then $S$ is the spectrum of the ring of integers of $K$; if $\Char K > 0$, then $S$ is the proper smooth curve over a finite field such that the function field of $S$ is $K$; a place of $K$ is denoted by $v$, and $K_v$ is the resulting completion; if $v$ is nonarchimedean, then $\cO_v$ and $\bF_v$ are the ring of integers and the residue field of $K_v$. We fix a set $\Sigma$ of places of $K$ that contains the archimedean places. For a finite subset $\Sigma_0 \subset \Sigma$ that contains the archimedean places, we set 
\[
\tst \bA_{K, \Sigma_0}^{\in \Sigma} \ce \prod_{v\in \Sigma_0} K_v \times \prod_{v\in \Sigma \setminus \Sigma_0} \cO_v,
\] 
to be interpreted as the zero ring if $\Sigma = \emptyset$. We let $\bA_K^{\in \Sigma}$ denote\footnote{We fear that the shorter notation $\bA_K^\Sigma$ is confusing: should $\bA_K^\Sigma$ stand for $\bA_K^{\in \Sigma}$ or for the similarly constructed~$\bA_K^{\not\in \Sigma}$?}
 the ring of ``$\Sigma$-adeles'', i.e., 
\[
\tst\bA_K^{\in \Sigma} \ce \varinjlim_{\Sigma_0 \subset \Sigma,\ \#\Sigma_0 < \infty} \bA_{K, \Sigma_0}^{\in \Sigma}.
\]
Thus, if $\Sigma$ is the set of all places, then $\bA_K^{\in \Sigma}$ is the adeles $\bA_K$, and if $\Sigma$ is finite, then $\bA_K^{\in\Sigma} = \prod_{v\in \Sigma} K_v$.

We identify the closed points of $S$ with the nonarchimedean $v$. For an open $U \subset S$, the notation $v\not\in U$ signifies that $v$ does not correspond to a point of $U$ (and hence could be archimedean). We denote by $U \cap \Sigma$ the set of those places in $\Sigma$ that correspond to a point of $U$. We write $S - \Sigma$ for the Dedekind scheme that is the localization of $S$ away from $\Sigma$, i.e., $S - \Sigma = \varprojlim U$ where $U$ ranges over the nonempty opens of $S$ with $S \setminus U \subset \Sigma$. We say that such a $U$ \emph{contains} $S - \Sigma$. 
\epp

\bpp[Conventions] \lab{conv}
For a field $F$, a choice of its algebraic closure $\ov{F}$ is made implicitly. Likewise implicit are the choices of embeddings $\ov{K} \hra \ov{K}_v$. All the cohomology other than profinite group cohomology is computed in the fppf topology, but we implicitly use \cite{Gro68c}*{11.7~$1^{\circ})$} (or \cite{Ces15d}*{B.17} in the algebraic space case) to make identifications with \'{e}tale cohomology when the coefficients are a smooth group scheme (or a smooth algebraic space); further identifications with Galois cohomology are likewise implicit. Derived functor $H^n$ is taken to vanish for $n < 0$. The neutral class of a nonabelian $H^1$ is denoted by $*$.

All actions and homogeneous spaces (and hence also torsors) are right instead of left unless noted otherwise. Algebraic spaces are not assumed to be quasi-separated. We implicitly use \cite{SP}*{\href{http://stacks.math.columbia.edu/tag/04SK}{04SK}}, i.e.,~the fppf local nature of being an algebraic space, to ensure that sheaves at hand, e.g.,~torsors under group algebraic spaces, are representable by algebraic spaces. For a scheme $T$, a $T$-group algebraic space $G$, and an fppf cover $T\pr/T$, the set of isomorphism classes of $G$-torsors trivialized by $T\pr$ is denoted by $H^1(T\pr/T, G)$; we implicitly use the cocycle point of view to identify $H^1(T\pr/T, G)$ with a subquotient of $G(T\pr\times_T T\pr)$ under an action of $G(T\pr)$.

In this paper `locally compact' means that every point has a compact neighborhood. Neither compact nor locally compact spaces are assumed to be Hausdorff. The Pontryagin dual of a locally compact Hausdorff abelian topological group $A$ is denoted by $A^*$ and is defined by $A^* \ce \Hom_\cont(A, \bR/\bZ)$ (we will always discuss Pontryagin duality in the context of torsion groups, or even $n$-torsion groups for a fixed $n \in \bZ_{\ge 1}$, so we could equivalently use $A^* \ce \Hom_\cont(A, \bQ/\bZ)$ as our definition). Topology on cohomology of local fields (resp.,~on cohomology of adele rings) is always the one defined in \cite{Ces15d}*{3.1--3.2} (resp.,~in \S\ref{ad-topo-Hn}).

\epp

\subsection*{Acknowledgements}
I thank Brendan Creutz for questions and suggestions that led to Poitou--Tate and Cassels--Poitou--Tate sequences proved in this paper. I thank the referees for helpful comments and suggestions. I thank James Milne for pointing out \Cref{Azu-rem} and for sending me a scan of \cite{Azu66}. I thank Brendan Creutz, Cristian Gonz{\'a}lez-Avil{\'e}s, Bjorn Poonen, and Yunqing Tang for helpful conversations or correspondence regarding the material of this paper. Some of the results presented here were obtained while the author was a Research Fellow at the Miller Institute for Basic Research in Science at the University of California Berkeley. I thank the Miller Institute and UC Berkeley for excellent conditions for~research.


\section{Cohomology of the adeles as a restricted product} \lab{coho-res}

The goal of this section is to prove in \Cref{res-pr-comp} that, under suitable hypotheses on the coefficients, the cohomology of the adeles $\bA_K^{\in \Sigma}$ decomposes as a restricted product of local cohomologies. Such a decomposition is most interesting when the coefficients are a $K$-group scheme $G$, for instance, in \S\S\ref{disc-image}--\ref{CPT} we will use it for a commutative finite $G$ (see also \Cref{no-vb,CFT} for examples of concrete conclusions with $G = \GL_n$ or $G = \bG_m$).

The proof begins by reducing to the case of a smooth $G$ and then proceeds by induction on the degree of the cohomology in question, with some of the key steps happening in the proof of \Cref{sm-key}. In these key steps we apply the inductive hypothesis to coefficients that need not come from $K$ by base change. More precisely, the relevant coefficients are quotients of restrictions of scalars along $R'/R$, where $R$ is an infinite product of rings of integers of nonarchimedean completions of $K$ and $R'$ is an infinite product of rings of integers of unramified extensions of those completions. 

In spite of not being base changes from $K$, these restrictions of scalars are nevertheless schemes due to the quasi-projectivity of $G$. The inductive argument therefore requires applying some of the intermediate results to a priori arbitrary $\bA_K^{\in \Sigma}$-group schemes $\cG$ of finite presentation. However, allowing such $\cG$ leads to further complications: due to the loss of quasi-projectivity, we no longer know how to argue the representability by schemes of the new restrictions of scalars and of their quotients. (Restricting to affine $\cG$ would resolve some of these difficulties but would exclude the possibility of $G$ being a semiabelian variety.) In contrast, the representability by algebraic spaces holds, but comes at the expense of also having to allow algebraic spaces in the intermediate results.

In conclusion, for purely technical reasons we formulate our results for coefficients that are $\bA_K^{\in\Sigma}$-group algebraic spaces. This generality makes the overall proof work, but does not lead to new cases of interest because quasi-separated $K$-group algebraic spaces are automatically schemes (as is recalled in \Cref{Art-input}).

To avoid cluttering the subsequent arguments with repetitive citations, we begin by recording the following two lemmas that concern limit arguments for algebraic spaces. We recall that for a scheme $T$ and a $T$-group algebraic space $\cG$, viewed as the sheaf of groups on the big fppf site of $T$, the notation $H^n(T, \cG)$ stands for $\cG(T)$ if $n = 0$, for the pointed set of isomorphism classes of $\cG$-torsor fppf sheaves (that are necessarily $T$-algebraic spaces) if $n = 1$, and for the derived functor fppf cohomology with coefficients $\cG$ if $n \ge 1$ and $\cG$ is commutative (in the commutative case, the two possible meanings of $H^1(T, \cG)$ agree by \cite{Gir71}*{III.3.5.4}).

\blem \lab{lim}
Let $I$ be a directed partially ordered set, $(T_i)_{i \in I}$ an inverse system of quasi-compact and quasi-separated schemes, $i_0 \in I$, and $\cG_{i_0}$ a $T_{i_0}$-group algebraic space that is locally of finite presentation. Suppose that the maps $T_{i\pr} \ra T_i$ are affine, set $T \ce \varprojlim T_i$ (which is a scheme by \cite{EGAIV3}*{8.2.3}), and set $\cG \ce (\cG_{i_0})_T$ and $\cG_i \ce (\cG_{i_0})_{T_i}$ for $i \ge i_0$. The pullback maps induce isomorphisms
\[
\textstyle \varinjlim_{i \ge i_0} H^n(T_i, \cG_i) \cong H^n(T, \cG)
\]
for $n \le 1$, and also for $n \ge 2$ if $\cG_{i_0}$ is commutative.
\elem

\bpf 
We treat the case $n = 0$, the case $n = 1$, and the case $n \ge 2$ in \ref{lim-0}, \ref{lim-1}, and \ref{lim-2}, respectively.
\benuma \addtocounter{enumi}{-1}
\item \lab{lim-0}
In the case of an affine $T_{i_0}$, the claim is a special case of the functor of points criterion \cite{SP}*{\href{http://stacks.math.columbia.edu/tag/04AK}{04AK}} for being locally of finite presentation. To reduce the general case to the affine one, we cover $T_{i_0}$ by finitely many affine opens, cover their pairwise intersections by finitely many further affine opens, use the description of $T$ provided by \cite{EGAIV3}*{8.2.2}, and use the fact that filtered direct limits commute with finite inverse limits in the category of sets.

\item \lab{lim-1}
For a faithfully flat morphism $T_{i_0}\pr \ra T_{i_0}$ of finite presentation and $i \ge i_0$, we set 
\[
\quad T\pr_i \ce (T_{i_0}\pr)_{T_i} \qq \text{and} \qq T\pr \ce (T_{i_0}\pr)_T.
\] 
By the $n = 0$ case applied to $(T_{i}\pr)_{i \in I}$, to $(T_{i}\pr\times_{T_i} T_i\pr)_{i \in I}$, and to $(T_{i}\pr\times_{T_i} T_i\pr\times_{T_i} T_i\pr)_{i \in I}$, 
\be \lab{Cech-H1}
\textstyle\varinjlim_{i \ge i_0} H^1(T_i\pr/T_i, \cG_i) \cong H^1(T\pr/T, \cG),
\ee
where we have again used the fact that taking filtered direct limits is exact. Since each $T_i$ (resp.,~$T$) is quasi-compact and every torsor under $\cG_i$ (resp.,~$\cG$) is an algebraic space that is faithfully flat and locally of finite presentation over $T_i$ (resp.,~$T$), every such torsor is trivialized by a single scheme morphism that is faithfully flat and of finite presentation. The claim therefore results by varying the $i_0$ and the $T\pr_{i_0}$ in \eqref{Cech-H1}.

\item \lab{lim-2}
The proof of \cite{SGA4II}*{VII, 5.9} based on topos-theoretic generalities of \cite{SGA4II}*{VI} and written for \'{e}tale cohomology with scheme coefficients continues to work for fppf cohomology with algebraic space coefficients (reasoning as in the $n = 0$ case for the limit formalism).
\qedhere
\eenum
\epf

\blem \lab{limit}
In the setup of \Cref{lim}, let $X_{i_0}$ be a $T_{i_0}$-algebraic space of finite presentation. Set 
\[
X \ce (X_{i_0})_{T} \qq \text{and} \qq X_{i} \ce (X_{i_0})_{T_i} \quad \text{for $i \ge i_0$.}
\]
If $X \ra T$ is separated (resp.,~is proper, resp.,~is flat, resp.,~is smooth, resp.,~has geometrically connected fibers), then the same holds for $X_{i} \ra T_i$ for large enough $i$.
\elem

\bpf
For all the properties except fibral geometric connectedness, the claim is a special case of \cite{Ryd15}*{B.3}. For fibral geometric connectedness, we apply \cite{Con12b}*{5.1} and \cite{EGAIV3}*{8.3.4}.
\epf

We now turn to the setup needed for the sought restricted product interpretations of cohomology.

\bpp[The comparison maps] \lab{app-setup}
Let $\cG$ be an $\bA_K^{\in \Sigma}$-group algebraic space of finite presentation. Limit arguments supplied by \cite{SP}*{\href{http://stacks.math.columbia.edu/tag/07SK}{07SK}} descend $\cG$ to a (commutative if so is $\cG$) $\bA_{K, \Sigma_0}^{\in \Sigma}$-group algebraic space $\cG_{\Sigma_0}$ of finite presentation for some finite subset $\Sigma_0 \subset \Sigma$ containing all the infinite places and prove that $\cG_{\Sigma_0}$ is unique up to base change to $\bA_{K, \Sigma\pr_0}^{\in \Sigma}$ with $ \Sigma_0 \subset \Sigma\pr_0 \subset \Sigma$ and $\#\Sigma\pr_0 < \infty$ and up to isomorphism. 

The main example is $\cG = G_{\bA_K^{\in \Sigma}}$ for a $K$-group scheme $G$ of finite type. For such $\cG$ we may choose a $U$-model of $G$ of finite presentation for a nonempty open $U \subset S$, let $\Sigma_0$ contain all the places of $\Sigma$ outside $U$, and let $\cG_{\Sigma_0}$ be the base change of the chosen $U$-model. 

By \Cref{lim}, for $n \le 1$, and also for $n \ge 2$ if $\cG$ is commutative, the pullback maps 
\be \lab{pr-pullb}
\textstyle H^n(\bA_{K, \Sigma\pr_0}^{\in \Sigma}, \cG_{\Sigma_0}) \ra \prod_{v \in \Sigma\pr_0} H^n(K_v, \cG) \times \prod_{v\in \Sigma \setminus \Sigma\pr_0} H^n(\cO_v, \cG_{\Sigma_0})
\ee
for varying $\Sigma_0\pr$ give rise to the comparison map
\be \lab{ad-pullb}
\textstyle H^n(\bA_K^{\in \Sigma}, \cG) \ra \varinjlim_{\Sigma_0 \subset \Sigma\pr_0 \subset \Sigma,\ \#\Sigma_0\pr < \infty} \p{\prod_{v \in \Sigma\pr_0} H^n(K_v, \cG) \times \prod_{v\in \Sigma\setminus \Sigma_0\pr} H^n(\cO_v, \cG_{\Sigma_0})}, 
\ee
whose target does not depend on the choices of $\Sigma_0$ and $\cG_{\Sigma_0}$.
\epp

We seek to prove in \Cref{ad-sm-main,ad-main} that the comparison map \eqref{ad-pullb} is bijective in many cases, in particular, if $\cG = G_{\bA_K^{\in\Sigma}}$. 
Our methods will exploit the following lemma, whose case of local $R_i$ would suffice for us.

\blem[\cite{Bha14}*{Thm.~1.3}] \lab{Bha-input}
For a set of rings $\{R_i\}_{i \in I}$ and a quasi-compact and quasi-separated $\bZ$-algebraic space $X$, the natural map 
\[
\textstyle X(\prod_i R_i) \ra \prod_i X(R_i)
\]
is bijective. \QED
\elem

The proofs in this section could be simplified significantly if \Cref{Bha-input} stayed true for $X$ that are algebraic stacks instead of algebraic spaces. However, \Cref{Bha-input} fails for such $X$: see \cite{Bha14}*{8.4}.

\Cref{Bha-input} allows us to settle some initial cases of the bijectivity of the comparison map.

\bprop \lab{H1-inj}
For $n = 0$ (resp.,~$n = 1$), the maps \eqref{pr-pullb}--\eqref{ad-pullb} are bijective (resp.,~injective). 
\eprop

\bpf
It suffices to treat \eqref{pr-pullb}. The $n = 0$ case is immediate from \Cref{Bha-input}. 

In the $n = 1$ case, a standard twisting argument \cite{Gir71}*{III.2.6.3} reduces us to proving that a torsor under an inner twist of $(\cG_{\Sigma_0})_{\bA_{K, \Sigma\pr_0}^{\in \Sigma}}$ is trivial as soon as its pullbacks to $K_v$ for $v\in \Sigma_0'$ and to $\cO_v$ for $v \in \Sigma \setminus \Sigma_0'$ are trivial. Such triviality follows by applying \Cref{Bha-input} because inner twists of $(\cG_{\Sigma_0})_{\bA_{K, \Sigma\pr_0}^{\in \Sigma}}$ and their torsors are algebraic spaces by \cite{SP}*{\href{http://stacks.math.columbia.edu/tag/04SK}{04SK}}, and they inherit quasi-compactness and quasi-separatedness by \cite{SP}*{\href{http://stacks.math.columbia.edu/tag/041L}{041L} and \href{http://stacks.math.columbia.edu/tag/041N}{041N}}.
\epf

\bcor \lab{no-vb}
There are no nontrivial vector bundles on $\Spec \bA_K$.
\ecor

\bpf
\Cref{H1-inj} allows us to infer the vanishing of $H^1(\bA_K, \GL_m)$ from the vanishing of $H^1(\cO_v, \GL_m)$ for $v\nmid \infty$ and of $H^1(K_v, \GL_m)$.
\epf

We start building up for \Cref{sm-key}, which is a key step towards proving the bijectivity of the comparison map for smooth $\cG$. \Cref{Lang-revamp}, whose part \ref{Lang-revamp-b} extends Lang's theorem \cite{Lan56}*{Thm.~2}, will be helpful throughout \S\S\ref{coho-res}--\ref{topo-ad}.

\blem \lab{Lang-revamp}
Let $k$ be a field, $G$ a finite type $k$-group scheme, and $X$ a homogeneous space~under~$G$.
\benum
\item \lab{Lang-revamp-a}
If $k$ is perfect and $G$ is finite connected, then $\#X(k) = 1$, so that, in particular, every $G$-torsor is trivial.

\item \lab{Lang-revamp-b}
If $k$ is finite and $X$ is geometrically connected, then $X(k) \neq \emptyset$. 

\item \lab{Lang-revamp-c}
If $k$ is finite and $G$ is connected, then $H^1(k, G) = \{ * \}$.

\item \lab{ff-H1-surj}
If $k$ is perfect and of dimension $\le 1$ (in the sense of \cite{Ser02}*{II.\S3.1}) and $G$ fits into an exact sequence 
\be \lab{Spr}
1 \ra G\pr \ra G \ra G^{\prime\prime} \ra 1
\ee
of $k$-group schemes of finite type, then the map 
\[
H^1(k, G) \ra H^1(k, G^{\prime\prime})
\]
is surjective.

\item \lab{Lang-revamp-d}
If $k$ is finite and $G$ is commutative, then $H^n(k, G) = 0$ for $n \ge 2$.
\eenum
\elem

\bpf \hfill
\benum
\item
We use the \'{e}tale sheaf property of $X$ to reduce to the $X(k) \neq \emptyset$ case by passing to a finite, and hence separable, extension $k\pr$ of $k$. Then we enlarge $k$ further to reduce to the $k = \ov{k}$ case. If $k = \ov{k}$, then the singleton $G(k)$ acts transitively on $X(k)$, and the claim follows.

\item
Due to geometric connectedness, $X$ is a homogeneous space under the identity component $G^0$, so we assume that $G$ is connected. Then, by \cite{SGA3Inew}*{\upshape{VII}$_{\text{\upshape{A}}}$, 8.3}, $G$ is an extension of a smooth connected $k$-group scheme $Q$ by a finite connected $H$. Since $X/H$ is a homogeneous space under $Q$, \cite{Ser02}*{III.\S2.4, Thm.~3} and Lang's theorem \cite{Lan56}*{Thm.~2} let us choose a $p \in (X/H)(k)$. The fiber $X_p$ is a homogeneous space under $H$, so $X_p(k)\neq \emptyset$ by \ref{Lang-revamp-a}.

\item
It suffices to apply \ref{Lang-revamp-b} to principal homogeneous spaces (i.e.,~torsors) under $G$.

\item
In the case when $G$ and $G^{\prime\prime}$ are smooth, the claim is a corollary of a theorem of Springer, see \cite{Ser02}*{III.\S2.4, Cor.~2}. We will reduce the general case to this smooth case.

By \cite{CGP10}*{C.4.1}, there is a $k$-smooth closed $k$-subgroup scheme $\wt{G^{\prime\prime}} \subset G^{\prime\prime}$ for which $\wt{G^{\prime\prime}}(L) = G^{\prime\prime}(L)$ for every finite (automatically separable) field extension $L$ of $k$. The description of torsors through cocycles therefore gives the induced bijection
\[
\qq H^1(k, \wt{G^{\prime\prime}}) \isomto H^1(k, G^{\prime\prime}).
\]
Thus, by pulling back \eqref{Spr} along $\wt{G^{\prime\prime}} \subset G^{\prime\prime}$, we reduce to the case when $G^{\prime\prime}$ is smooth.

To reduce further to the case when $G'$ is also smooth, we begin by applying \cite{SGA3Inew}*{VII$_{\text{\upshape{A}}}$,~8.3} to infer the existence of a finite connected subgroup $N \subset G'$ that is normal even in $G$ and for which $G'/N$ is smooth ($N$ is either trivial or a large enough Frobenius kernel). Then, by \ref{Lang-revamp-a}, 
\[
\qqq G(L) = (G/N)(L)
\]
for every finite field extension $L$ of $k$, so we again get a bijection
\[
\qqq H^1(k, G) \isomto H^1(k, G/N).
\]
Passage to quotients by $N$ therefore reduces us to the case when both $G'$ and $G^{\prime\prime}$ are smooth. Then, by \cite{SGA3Inew}*{9.2~(viii)}, $G$ is also smooth, which completes the promised reduction.

\item
If $G$ is smooth, then a cohomological dimension argument suffices. In general, \cite{Ces15d}*{6.4} embeds $G$ into a commutative smooth $\wt{G}$, and \ref{ff-H1-surj} allows replacing $G$ by the smooth~$\wt{G}/G$.
\qedhere
\eenum
\epf

We will sometimes use the following lemma implicitly. 

\blem[\cite{Art69b}*{4.2}] \lab{Art-input}
For a field $k$, every quasi-separated $k$-group algebraic space locally of finite type is a scheme. \QED
\elem

\bcor \lab{Hens-0}
For a Henselian local ring $R$ with a finite residue field $\bF$ and a flat $R$-group algebraic space $G$ of finite presentation, $H^n(R, G)$ vanishes if $n = 1$ and $G$ is smooth with connected $G_{\bF}$, and also if $n \ge 2$ and $G$ is commutative.
\ecor

\bpf
By \cite{Ces15d}*{B.11} and \cite{Ces15d}*{B.13} (which is based on \cite{Toe11}*{3.4}), 
\[
H^n(R, G) \ra H^n(\bF, G)
\]
is injective, so \Cref{Lang-revamp} \ref{Lang-revamp-c} and \ref{Lang-revamp-d} (together with \Cref{Art-input}) give the claim.
\epf

\blem \lab{res-qcqs}
Let $T\pr \ra T$ be a finite \'{e}tale morphism of schemes. For a $T\pr$-algebraic space $X\pr$, let $X$ denote the restriction of scalars $\Res_{T\pr/T} X\pr$ (which is a $T$-algebraic space by \cite{SP}*{\href{http://stacks.math.columbia.edu/tag/05YF}{05YF}}). If a $T\pr$-morphism $f\pr\colon X\pr \ra Y\pr$ is quasi-compact (resp.,~quasi-separated), then so is $\Res_{T\pr/T} f\pr \colon X \ra Y$.
\elem

\bpf
Since $\Res_{T\pr/T}$ commutes with fiber products and base change in $T$, we only need to prove the quasi-compactness claim and may work \'{e}tale locally on $T$. Limit arguments therefore permit us to assume that $T\pr = \bigsqcup_{j = 1}^d T$, in which case $\Res_{T\pr/T} f\pr$ is a product of $d$ quasi-compact $T$-morphisms.
\epf

\blem \lab{construct}
Let $T$ be a quasi-compact scheme and $X$ a $T$-algebraic space of finite presentation. There is an $m_X \in \bZ_{>0}$ such that for each geometric fiber $X_{\ov{t}}$ of $X \ra T$, the underlying topological space $\abs{X_{\ov{t}}}$ has at most $m_X$ connected components.
\elem

\bpf
We replace $T$ by an affine open to reduce to the case when $T$ is affine. In this case, we take an \'{e}tale cover $Y \surjects X$ of finite presentation with $Y$ a scheme and replace $X$ by $Y$ to reduce to the scheme case, which is settled by \cite{EGAIV3}*{9.7.9} (see also \cite{EGAIII1}*{\b{0}.9.3.1}).
\epf

We are now ready to use the preceding auxiliary lemmas to settle a key step towards the bijectivity of the comparison map.

\blem \lab{sm-key}
Fix a set $I$. For $i \in I$, let $R_i$ be a Henselian local ring with finite residue field $\bF_i$. For a smooth $\p{\prod_{i} R_i}$-group algebraic space $\cF$ of finite presentation, the pullback map
\[
\textstyle H^n(\prod_{i} R_i, \cF) \ra \prod_i H^n(R_i, \cF)
\]
is bijective for $n\le 1$, and also for $n \ge 2$ if $\cF$ is commutative. For $n \ge 2$, both 
\[
\tst H^n(\prod_{i} R_i, \cF) = 0 \qq \text{and} \qq \prod_i H^n(R_i, \cF) = 0;
\]
the same holds for $n = 1$ if each $\cF_{\bF_i}$ is connected.
\elem

\bpf
The last sentence follows from the rest due to \Cref{Hens-0}.

For a fixed $d \in \bZ_{> 0}$, let $R_i\pr$ be the Henselian local finite \'{e}tale $R_i$-algebra that corresponds to the residue field extension $\bF_i\pr/\bF_i$ of degree $d$. Set 
\[
\tst R \ce \prod_{i} R_i \qq \text{and} \qq R\pr \ce \prod_i R_i\pr,
\]
so that $R\pr/R$ is finite \'{e}tale of degree $d$ by \cite{Con12a}*{7.5.5}. We consider the restriction of scalars 
\[
\cH \ce \Res_{R\pr/R}(\cF_{R\pr}),
\]
which is an $R$-group algebraic space by \cite{SP}*{\href{http://stacks.math.columbia.edu/tag/05YF}{05YF}}. The adjunction map $\iota\colon \cF \ra \cH$ is injective, so 
\[
\cQ \ce \cH/\iota(\cF)
\] 
is an $R$-algebraic space by \cite{SP}*{\href{http://stacks.math.columbia.edu/tag/06PH}{06PH}}.

\bcl \lab{HQ-qcqs}
Both $\cH$ and $\cQ$ are smooth and of finite presentation over $R$.
\ecl

\bpf
Smoothness of $\cH$ follows from \cite{SP}*{\href{http://stacks.math.columbia.edu/tag/04AK}{04AK} and \href{http://stacks.math.columbia.edu/tag/04AM}{04AM}} and implies that of $\cQ$ by \cite{SP}*{\href{http://stacks.math.columbia.edu/tag/0AHE}{0AHE}}. Thus, $\cH$ is of finite presentation by \Cref{res-qcqs}, and hence so is $\cQ$ by \cite{SP}*{\href{http://stacks.math.columbia.edu/tag/040W}{040W} and \href{http://stacks.math.columbia.edu/tag/03KS}{03KS}}. 
\epf

\bcl \lab{H1-sm-surj}
The map $\textstyle H^1(\prod_{i} R_i, \cF) \ra \prod_{i} H^1(R_i, \cF)$ is surjective.
\ecl

\bpf
We choose $m_\cF$ as in \Cref{construct} (applied with $X = \cF$) and set $d \ce m_\cF!$ in the construction of $\cH$ and $\cQ$. Thanks to Lang's theorem \cite{Lan56}*{Thm.~2}, $H^1(\bF_i, \cF)$ and $H^1(\bF_i\pr, \cF)$ inject into analogous $H^1$'s with the component group of $\cF_{\bF_i}$ as coefficients. Every torsor under this component group is a finite \'{e}tale $\bF_i$-scheme of degree at most $m_\cF$, so it has an $\bF_i\pr$-point. In conclusion,
\[
\tst H^1(\bF_i, \cF) \ra H^1(\bF_i\pr, \cF)
\] 
vanishes, and hence so does each 
\[
\tst H^1(R_i, \cF) \ra H^1(R_i\pr, \cF)
\]
by \cite{Ces15d}*{B.11}. Therefore, \cite{Gir71}*{V.3.1.3} gives the vanishing of $b$ in
\be \lab{SES}
\cQ(R_i) \xra{\delta_i} H^1(R_i, \cF) \xra{b} H^1(R_i, \cH),
\ee
and the surjectivity of $\delta_i$ follows. It remains to use \Cref{HQ-qcqs} to apply \Cref{Bha-input} to $\cQ$ and use the compatibility of the sequences \eqref{SES} with their counterpart over $\prod_{i} R_i$.
\epf

As in \Cref{H1-inj}, \Cref{Bha-input} and \cite{Gir71}*{III.2.6.3} prove that the map
\[
\tst H^n(\prod_{i} R_i, \cF) \ra \prod_{i} H^n(R_i, \cF)
\]
is bijective for $n = 0$ and injective for $n = 1$, so \Cref{H1-sm-surj} settles the $n \le 1$~case. 

Suppose from now on that $n \ge 2$ and $\cF$ is commutative, so $\prod_i H^n(R_i, \cF) = 0$ by \Cref{Hens-0}. We argue the remaining vanishing of $H^n(R, \cF)$ by induction on $n$ (using the settled bijectivity in the $n = 1$ case to get started).

We fix an $\gA \in H^n(R, \cF)$ and seek to prove that $\gA = 0$. \Cref{lim} applied to strict Henselizations of $\Spec R$ provides an \'{e}tale cover 
\[
p\colon \Spec A \surjects \Spec R \qq  \text{with $\gA|_A = 0$.}
\]
By \Cref{construct}, the fibral degrees of $p$ are bounded by some $m \in \bZ_{> 0}$. We set $d\ce m!$, and, as before, let $R_i'$ be the finite \'{e}tale $R_i$ algebra corresponding to the residue field extension of degree $d$ and set $R' \ce \prod_i R_i'$. Each $\Spec R_i' \ra \Spec R$ factors through $p$, and hence so does $p\pr\colon \Spec R\pr \ra \Spec R$, so that $\gA|_{R\pr} = 0$. The acyclicity of $p\pr_*$ supplied by \cite{SGA4II}*{VIII, 5.5} shows that the spectral sequence 
\[
H^n(R, \bbR^m p'_*( \cF_{R'})) \Rightarrow H^{n + m}(R', \cF)
\]
degenerates, so the equality $\gA|_{R'} = 0$ proves that the map
\[
\tst H^n(R, \cF) \ra H^n(R, \cH) \cong H^n(R', \cF)
\] 
kills $\gA$. Finally, \Cref{HQ-qcqs} and the inductive hypothesis provide the isomorphisms in
\[
\xymatrix{
H^{n - 1}(R, \cH) \ar[r] \ar[d]^-{\cong} & H^{n - 1}(R, \cQ) \ar[r] \ar[d]^-{\cong}  & H^n(R, \cF) \ar[r]\ar[d] & H^n(R, \cH) \ar[d] \\
\prod_i H^{n - 1}(R_i, \cH) \ar[r] & \prod_i H^{n - 1}(R_i, \cQ) \ar[r]  & 0 \ar[r] & 0,
}
\]
so a diagram chase proves that $\gA = 0$.
\epf

With \Cref{sm-key} at hand, we are ready to prove the bijectivity of the comparison map for smooth~$\cG$.

\bthm \lab{ad-sm-main}
For a smooth $\bA_K^{\in \Sigma}$-group algebraic space $\cG$ of finite presentation, the comparison~map 
\be\lab{ad-sm-eq}
\textstyle H^n(\bA_K^{\in \Sigma}, \cG) \ra \varinjlim_{\Sigma_0\pr} \p{ \prod_{v \in \Sigma\pr_0} H^n(K_v, \cG) \times \prod_{v\in \Sigma \setminus \Sigma_0\pr} H^n(\cO_v, \cG_{\Sigma_0})}
\ee
of \eqref{ad-pullb} is bijective for $n \le 1$, and also for $n \ge 2$ if $\cG$ is commutative. For $n \ge 2$, this reduces to 
\[
\tst H^n(\bA_K^{\in\Sigma}, \cG) \cong \bigoplus_{v\in\Sigma} H^n(K_v, \cG);
\]
the same holds for $n = 1$ if $\cG \ra \Spec \bA_K^{\in\Sigma}$ has connected fibers.
\ethm

\bpf
We apply \Cref{limit} to ensure that $\cG_{\Sigma_0} \ra \Spec \bA_{K, \Sigma_0}^{\in\Sigma}$ of \S\ref{app-setup} inherits smoothness or fibral connectedness from $\cG \ra \Spec \bA_K^{\in\Sigma}$. For a fixed $\Sigma\pr_0$, \Cref{sm-key} and the decomposition 
\[
\textstyle\Spec \bA_{K, \Sigma\pr_0}^{\in\Sigma} = \p{\bigsqcup_{v\in \Sigma\pr_0} \Spec K_v} \sqcup \Spec(\prod_{v\in\Sigma \setminus\Sigma\pr_0} \cO_v)
\] 
show that \eqref{pr-pullb} is bijective, so the bijectivity of \eqref{ad-pullb} (i.e.,~of \eqref{ad-sm-eq}) follows. For the last claim, we apply the vanishing aspect of \Cref{sm-key}.
\epf

The following corollary, mentioned to us by Bjorn Poonen, was one of motivations to work out the general \Cref{ad-sm-main,ad-main}. We recall that the Brauer group $\Br(T)$ of a scheme $T$ is by definition the group of similarity classes of Azumaya algebras over $T$.

\bcor \lab{CFT}
The Albert--Brauer--Hasse--Noether exact sequence from global class field theory may be rewritten as
\[
0 \ra \Br(K) \ra \Br(\bA_K) \ra \bQ/\bZ \ra 0.
\]
\ecor

\bpf
\Cref{ad-sm-main} with $\cG = \bG_m$ gives 
\[
\tst H^2(\bA_K, \bG_m) \cong \bigoplus_v \Br(K_v), \qq \text{so that} \qq H^2(\bA_K, \bG_m) = H^2(\bA_K, \bG_m)_\tors,
\] 
which allows us to apply \cite{Gab81}*{Ch.~II, Thm.~1 on p.~163} to conclude that also 
\[
\tst \Br(\bA_K) \cong H^2(\bA_K, \bG_m). \qedhere
\]
\epf

\brem \lab{Azu-rem}
In \cite{Azu66}, assuming that $K$ is a number field, Azumaya gave a direct computation of $\Br(\bA_K)$ via central simple algebras and also proved that $\Br(\bA_K) \cong \bigoplus_v \Br(K_v)$.
\erem

We turn to the remaining input to our main result (\Cref{ad-main}) regarding the comparison map.

\blem \lab{uniform}
For a $K$-group scheme $G$ of finite type, there is a $d \in \bZ_{> 0}$, a nonempty open $U \subset S$, and a model $\cG \ra U$ of $G \ra \Spec K$ such that for every $v \in U$ there is an extension $K_v\pr/K_v$ of degree $d$ for which the pullback map
\[
\tst H^1(\cO_v, \cG) \ra H^1(\cO_v\pr, \cG)
\]
vanishes (here $\cO_v\pr$ is the ring of integers of $K_v\pr$).
\elem

\bpf
We begin with some preliminary observations.

\begin{enumerate}[label=(\arabic*)]
\item \lab{u-1}
For a nonarchimedean local field $k$, its ring of integers $\fo$, and a finite \'{e}tale $\fo$-group scheme $\cF$ of order $\#\cF$, there is an extension $k\pr/k$ of degree $(\#\cF)!$ for which the pullback map 
\[
\quad H^1(\fo, \cF) \ra H^1(\fo\pr, \cF)
\] 
vanishes (where $\fo\pr$ is the ring of integers of $k\pr$).
\eenum

\bpf
Every torsor under $\cF$ is represented by a finite \'{e}tale $\fo$-group scheme of order $\#\cF$. Thus, if we let $k'$ be the unramified extension of $k$ of degree $(\#\cF)!$, then every such torsor has an $\fo'$-point.
\epf

\benuma \addtocounter{enumi}{1}
\item \lab{u-3}
If a finite $\cF$ in \ref{u-1} is not assumed to be \'{e}tale but instead $k$ is of characteristic $p > 0$ and $\cF_k$ is connected, then the same conclusion holds with $(\# \cF)!$ replaced by $p^{(\#\cF)!}$.
\eenum

\bpf
By the valuative criterion of properness, 
\[
H^1(\fo, \cF) \ra H^1(k, \cF)
\] 
is injective. Thus, the claim results from the vanishing of 
\[
H^1(k, \cF) \ra H^1(k^{1/p^{(\#\cF)!}}, \cF)
\]
supplied by \cite{Ces15d}*{5.7 (b)}. 
\epf

\benuma \addtocounter{enumi}{2}
\item \lab{u-2}
For $k$ and $\fo$ as in \ref{u-1} and $\cC$ a smooth $\fo$-group scheme with connected fibers, 
\[
H^1(\fo, \cC) = \{*\}.
\]
\eenum

\bpf
This is a special case of \Cref{Hens-0}.
\epf

By \cite{SGA3Inew}*{\upshape{VI}$_{\text{\upshape{B}}}$, 1.6.1 et VII$_{\text{\upshape{A}}}$, 8.3}, $G$ is an extension of a smooth $K$-group $\wt{G}$ by a finite connected $K$-group $N$. In turn, $\wt{G}$ is an extension of its finite \'{e}tale component group by the smooth identity component $\wt{G}^0$. A short exact sequence of finite type $K$-groups spreads out, as does smoothness, finiteness, \'{e}taleness, or fibral connectedness of its terms (the latter by \cite{EGAIV3}*{9.7.8}), so the nonabelian cohomology exact sequences combine with \ref{u-1}--\ref{u-2} to prove the claim.
\epf

We are now ready to address the bijectivity of the comparison map in the $n = 1$ case.

\bprop \lab{H1-bij}
For $n = 1$ and $\cG \simeq G_{\bA_K^{\in\Sigma}}$ with $G$ a $K$-group scheme of finite type, the comparison map \eqref{ad-pullb} is bijective.
\eprop

\bpf
Let $d \in \bZ_{> 0}$ and a model $\cG \ra U$ of $G \ra \Spec K$ be as provided by \Cref{uniform}. We shrink $U$ to assume that $\cG \ra U$ is of finite presentation. For a $\Sigma_0$ containing all the places of $\Sigma$ outside $U$, we set $\cG_{\Sigma_0} \ce \cG_{\bA_{K, \Sigma_0}^{\in \Sigma}}$. Due to \Cref{H1-inj}, the surjectivity of \eqref{pr-pullb} argued below will suffice. 

To simplify the notation we may and do assume that $\Sigma\pr_0 = \Sigma_0$ in \eqref{pr-pullb}. Thanks to the decomposition 
\[
\tst \Spec(\bA_{K, \Sigma_0}^{\in \Sigma}) = \p{\bigsqcup_{v \in \Sigma_0} \Spec K_v} \sqcup \Spec (\prod_{v\in\Sigma \setminus \Sigma_0} \cO_v),
\]
we may focus on the surjectivity~of
\be\lab{surj-wanted}
\textstyle H^1(\prod_{v\in\Sigma \setminus \Sigma_0} \cO_v, \cG) \ra \prod_{v\in \Sigma \setminus \Sigma_0} H^1(\cO_v, \cG).
\ee
For $v\in \Sigma\setminus \Sigma_0$, we let $\cO_v\pr/\cO_v$ be a finite free extension of degree $d$ provided by \Cref{uniform}. We choose an $\cO_v$-basis of $\cO_v\pr$ for each $v\in \Sigma\setminus \Sigma_0$ to see that $\prod_{v\in \Sigma\setminus \Sigma_0} \cO_v\pr / \prod_{v\in \Sigma\setminus \Sigma_0} \cO_v$ is a finite free extension of degree $d$ and that the following natural homomorphisms are isomorphisms:
\[
\ba
\textstyle \prod_{v\in \Sigma\setminus \Sigma_0} \cO_v\pr \tensor_{\prod_{v\in \Sigma\setminus \Sigma_0} \cO_v} \prod_{v\in \Sigma\setminus \Sigma_0} \cO_v\pr &\isomto \textstyle\prod_{v\in \Sigma\setminus \Sigma_0} \cO_v\pr \tensor_{\cO_v} \cO_v\pr, \\
\textstyle\prod_{v\in \Sigma\setminus \Sigma_0} \cO_v\pr \tensor_{\prod_{v\in \Sigma\setminus \Sigma_0} \cO_v}\prod_{v\in \Sigma\setminus \Sigma_0} \cO_v\pr \tensor_{\prod_{v\in \Sigma\setminus \Sigma_0} \cO_v} \prod_{v\in \Sigma\setminus \Sigma_0} \cO_v\pr &\isomto \textstyle\prod_{v\in \Sigma\setminus \Sigma_0} \cO_v\pr \tensor_{\cO_v} \cO_v\pr \tensor_{\cO_v} \cO_v\pr.
\ea
\]
For $v\in \Sigma\setminus \Sigma_0$, let $\cT_v \ra \Spec \cO_v$ be a $\cG_{\cO_v}$-torsor. By the choice of $\cO_v\pr$, there is a descent datum 
\[
\tst g_v \in \cG(\cO_v\pr \tensor_{\cO_v} \cO_v\pr)
\] 
that gives rise to $\cT_v$. \Cref{Bha-input} and the displayed bijections assemble the $g_v$ to a descent datum 
\[
\tst g \in \cG(\prod_{v\in \Sigma\setminus \Sigma_0} \cO_v\pr \tensor_{\prod_{v\in \Sigma\setminus \Sigma_0} \cO_v} \prod_{v\in \Sigma\setminus \Sigma_0} \cO_v\pr).
\] 
By construction, the class of the $\cG$-torsor that results from $g$ maps under \eqref{surj-wanted} to the tuple whose $v$-component is the class of $\cT_v$. In conclusion, \eqref{surj-wanted} is surjective, as desired.
\epf

With all ingredients in place, we now turn to the promised bijectivity of the comparison map.

\bthm \lab{ad-main}
For a $K$-group scheme $G$ of finite type, the comparison map 
\[
\textstyle H^n(\bA_K^{\in\Sigma}, G) \ra \varinjlim_{\Sigma_0\pr} \p{ \prod_{v \in \Sigma_0\pr} H^n(K_v, G) \times \prod_{v\in \Sigma\setminus \Sigma_0\pr} H^n(\cO_v, \cG_{\Sigma_0}) }
\]
of \eqref{ad-pullb} is bijective for $n \le 1$, and also for $n \ge 2$ if $G$ is commutative. For $n \ge 2$, this reduces to 
\[
\tst H^n(\bA_K^{\in\Sigma}, G) \cong \bigoplus_{v\in\Sigma} H^n(K_v, G);
\] 
the same holds for $n = 1$ if $G$ is connected and smooth.
\ethm

\bpf
\Cref{H1-inj,H1-bij} supply the bijectivity in the $n \le 1$ case, so we assume that $n \ge 2$.

By \cite{Ces15d}*{6.4}, $G$ is a closed subgroup of a $K$-group scheme $H$ that is of finite type, smooth, and commutative. By \cite{SGA3Inew}*{\upshape{VI}$_{\text{\upshape{A}}}$, 3.2}, $H/G$ is a $K$-group scheme that inherits these properties of $H$. The exact sequence 
\[
0 \ra G \ra H \ra H/G \ra 0
\] 
spreads out to an exact sequence of finite type models over a nonempty open $U \subset S$, so the five lemma allows us to deduce the claimed bijectivity for $G$ from the bijectivity in the smooth case established in \Cref{ad-sm-main}.

For the last sentence, \Cref{ad-sm-main} suffices in the $n = 1$ case. In the $n \ge 2$ case, we shrink $U$ to assume that the $U$-models of $H$ and $H/G$ are smooth. We then combine \cite{Gro68c}*{11.7~$2^{\circ})$} with the five lemma to get 
\[
\qq H^n(\cO_v, \cG) \cong H^n(\bF_v, \cG) \qq \text{for $v \in U$,} 
\]
where $\cG$ is the chosen $U$-model of $G$. It remains to note that $H^n(\bF_v, \cG) = 0$ by \Cref{Lang-revamp} \ref{Lang-revamp-d}.
\epf

In the following lemma we record conditions under which the target of the comparison map is a restricted product.

\blem \lab{res-pr}
Let $\cG$ be a finitely presented $\bA_K^{\in \Sigma}$-group algebraic space, $n\in \bZ_{\ge 0}$, and assume that
\benuma
\item \lab{res-pr-1}
$n = 0$, and $\cG$ is separated; or

\item \lab{res-pr-4}
$n = 1$, and $\cG$ is either smooth with connected fibers or proper; or

\item \lab{res-pr-3}
$n \ge 2$, and $\cG$ is commutative.
\eenum
Use \Cref{limit} to enlarge a $\Sigma_0$ of \S\ref{app-setup} to ensure that $\cG_{\Sigma_0} \ra \Spec \bA_{K, \Sigma_0}^{\in \Sigma}$ inherits the assumed properties of $\cG \ra \Spec \bA_K^{\in\Sigma}$. Then the pullback map 
\[
H^n(\cO_v, \cG_{\Sigma_0}) \ra H^n(K_v, \cG)
\]
is injective for all $v\in \Sigma \setminus \Sigma_0$, the resulting restricted product $\prod\pr_{v\in \Sigma} H^n(K_v, \cG)$ does not depend on the choices of $\Sigma_0$ and $\cG_{\Sigma_0}$, and $\prod\pr_{v\in \Sigma} H^n(K_v, \cG)$ identifies with the target of the comparison map \eqref{ad-pullb}.
\elem

\bpf
The injectivity follows from the valuative criterion of separatedness \cite{SP}*{\href{http://stacks.math.columbia.edu/tag/03KU}{03KU}}, \Cref{Hens-0}, and the valuative criterion of properness \cite{SP}*{\href{http://stacks.math.columbia.edu/tag/0A40}{0A40}} applied to the algebraic space parametrizing isomorphisms between two fixed $(\cG_{\Sigma_0})_{\cO_v}$-torsors for $v\in \Sigma \setminus \Sigma_0$ (the fppf sheaf parametrizing torsor isomorphisms is fppf locally isomorphic to $(\cG_{\Sigma_0})_{\cO_v}$ itself, so it is representable by an algebraic space due to \cite{SP}*{\href{http://stacks.math.columbia.edu/tag/04SK}{04SK}}). The independence follows from the second sentence of \S\ref{app-setup}. The last claim follows from the rest.
\epf

We conclude \S\ref{coho-res} by combining the previous results to get the promised restricted product interpretation of the cohomology of $\bA_K^{\in\Sigma}$.

\bthm \lab{res-pr-comp}
Let $\cG$ be a finitely presented $\bA_K^{\in \Sigma}$-group algebraic space, $n\in \bZ_{\ge 0}$, and assume that 
\benuma
\item \lab{res-pr-comp-1}
$n = 0$, and $\cG$ is separated; or

\item \lab{res-pr-comp-4}
$n = 1$, and $\cG$ satisfies any of the following: is isomorphic to $G_{\bA_K^{\in\Sigma}}$ for a proper $K$-group scheme $G$, is smooth with connected fibers, or is smooth and proper; or

\item \lab{res-pr-comp-3}
$n \ge 2$, and $\cG$ is commutative and either smooth or isomorphic to $G_{\bA_K^{\in \Sigma}}$ for a $K$-group scheme~$G$.
\eenum
Then the comparison map \eqref{ad-pullb} reduces to a bijection towards the restricted product of \Cref{res-pr}:
\[
\textstyle H^n(\bA_K^{\in\Sigma}, \cG) \isomto \prod\pr_{v\in \Sigma} H^n(K_v, \cG).
\]
\ethm

\bpf
The claim follows by combining \Cref{res-pr} with \Cref{H1-inj} and with \Cref{ad-sm-main,ad-main}.
\epf


\section{Topology on cohomology of the adeles} \lab{topo-ad}

The results of \S\ref{coho-res} are useful in the study of a natural topology carried by cohomology of the adeles. This topology plays an important role in the arithmetic duality proofs of \S\S\ref{disc-image}--\ref{CPT}, so in \S\ref{topo-ad} we gather the topological input needed for those proofs. For completeness, we also include several results that are not used in this paper but that seem indispensable to a robust topological theory. 

We formulate the results in the setting where the coefficients are group algebraic spaces because this does not lead to any complications in the proofs and is convenient for \S\ref{seq-cont} and \Cref{open-maps}, where we would otherwise need to impose additional assumptions to guarantee the representability of $\cQ$ by a scheme. However, if desired, the reader may restrict to the case of schemes throughout.

\bpp[Topology on $H^n(\bA_K^{\in\Sigma}, \cG)$] \lab{ad-topo-Hn}
We suppose that $\cG$ of \S\ref{app-setup} is $\bA_K^{\in\Sigma}$-flat (e.g.,~that $\cG \simeq G_{\bA_K^{\in\Sigma}}$) and use \Cref{limit} to ensure the $\bA_{K, \Sigma_0}^{\in\Sigma}$-flatness of $\cG_{\Sigma_0}$. 

For $n \le 1$, and also for $n \ge 2$ if $\cG$ is commutative, each $H^n(K_v, \cG)$ with $v \in \Sigma$ and each $H^n(\cO_v, \cG_{\Sigma_0})$ with $v\in \Sigma\setminus \Sigma_0$ carries a ``$v$-adic'' topology defined as in \cite{Ces15d}*{3.1--3.2}. Each pullback map
\be \lab{pullb-open}
H^n(\cO_v, \cG_{\Sigma_0}) \ra H^n(K_v, \cG) \quad \quad \text{ with $v  \in \Sigma \setminus \Sigma_0$} 
\ee
is continuous and open by \cite{Ces15d}*{3.5 (c), 2.7 (viii), and 2.9 (e) (with A.2 (b))} (with \Cref{Art-input} for the required separatedness of $(\cG_{\Sigma_0})_{\bF_v}$). Thus, 
\[
\tst\varinjlim_{\Sigma\pr_0} \p{ \prod_{v\in \Sigma\pr_0} H^n(K_v, \cG) \times \prod_{v\in \Sigma \setminus \Sigma\pr_0} H^n(\cO_v, \cG_{\Sigma_0}) },
\] 
which carries the direct limit topology of product topologies, has open transition maps, and the~maps
\be\lab{open-tr}
\prod_{v\in \Sigma\pr_0} H^n(K_v, \cG) \times \prod_{v\in \Sigma \setminus \Sigma\pr_0} H^n(\cO_v, \cG_{\Sigma_0}) \ra \varinjlim \textstyle\p{ \prod_{v\in \Sigma\pr_0} H^n(K_v, \cG) \times \prod_{v\in \Sigma \setminus \Sigma\pr_0} H^n(\cO_v, \cG_{\Sigma_0}) }
\ee
are open, too. We endow $H^n(\bA_K^{\in\Sigma}, \cG)$ with the preimage topology via the comparison map
\be\lab{ad-topo-def}
\textstyle H^n(\bA_K^{\in\Sigma}, \cG) \ra \varinjlim_{\Sigma_0 \subset \Sigma\pr_0 \subset \Sigma,\, \#\Sigma_0\pr < \infty} \p{ \prod_{v\in \Sigma\pr_0} H^n(K_v, \cG) \times \prod_{v\in \Sigma \setminus \Sigma\pr_0} H^n(\cO_v, \cG_{\Sigma_0}) }
\ee
of \eqref{ad-pullb}. By the second sentence of \S\ref{app-setup}, the topology on the target of \eqref{ad-topo-def} does not depend on the choices of $\Sigma_0$ and $\cG_{\Sigma_0}$, and hence neither does the topology on $H^n(\bA_K^{\in\Sigma}, \cG)$.
\epp

\brems
\remi
See \Cref{ad-sm-main,ad-main} for some settings in which the map \eqref{ad-topo-def} is bijective.

\remi
For $n = 0$ and $\cG = G_{\bA_K}$ with $G$ a $K$-group scheme of finite type, the topology of \S\ref{ad-topo-Hn} recovers the usual adelic topology on $G(\bA_K)$ discussed, for instance, in \cite{Con12b}*{\S4}. 
\erems

We record conditions, under which the target of \eqref{ad-topo-def} is homeomorphic to a restricted product.

\bprop \lab{res-pr-top}
Let $\cG$ be a flat $\bA_K^{\in\Sigma}$-group algebraic space of finite presentation and $n \in \bZ_{\ge 0}$.
\benum
\item \lab{res-pr-top-a}
If either of the conditions \ref{res-pr-1}--\ref{res-pr-3} of \Cref{res-pr} holds, then the restricted product topology on $\prod_{v\in\Sigma}\pr H^n(K_v, \cG)$ does not depend on the choice of $\Sigma_0$ and $\cG_{\Sigma_0}$ and the resulting identification 
\[
\tst\qqq \prod_{v\in\Sigma}\pr H^n(K_v, \cG) \cong \varinjlim_{\Sigma_0 \subset \Sigma\pr_0 \subset \Sigma,\, \#\Sigma_0\pr < \infty} \p{ \prod_{v\in \Sigma\pr_0} H^n(K_v, \cG) \times \prod_{v\in \Sigma \setminus \Sigma\pr_0} H^n(\cO_v, \cG_{\Sigma_0}) }
\]
is a homeomorphism.

\item \lab{res-pr-top-b}
If either of the conditions \ref{res-pr-comp-1}--\ref{res-pr-comp-3} of \Cref{res-pr-comp} holds, then the resulting identification 
\[
\tst \quad H^n(\bA_K^{\in\Sigma}, \cG) \cong \prod\pr_{v\in \Sigma} H^n(K_v, \cG)
\]
is a homeomorphism.
\eenum
\eprop

\bpf
The claims follow by combining the openness of the maps \eqref{pullb-open}--\eqref{open-tr} with their injectivity supplied by \Cref{res-pr}.
\epf

In \Cref{basic,twist,disc,haus,lc}, we record several basic properties of the adelic topology of \S\ref{ad-topo-Hn}.

\bprop \lab{basic}
Let $\cG$ and $\cG\pr$ be flat $\bA_K^{\in\Sigma}$-group algebraic spaces of finite presentation.
\benum
\item \lab{basic-a}
For $n \le 1$, and also for $n \ge 2$ if $\cG$ and $\cG\pr$ are commutative, the map 
\[
\tst H^n(\bA_K^{\in\Sigma}, \cG) \ra H^n(\bA_K^{\in\Sigma}, \cG\pr)
\] 
induced by an $\bA_K^{\in\Sigma}$-homomorphism $\cG \ra \cG\pr$ is continuous.

\item \lab{basic-b}
For $n \le 1$, and also for $n \ge 2$ and $\cG$ and $\cG\pr$ are commutative, the maps of \ref{basic-a} induce a homeomorphism 
\[
H^n(\bA_K^{\in\Sigma}, \cG \times \cG\pr) \isomto H^n(\bA_K^{\in\Sigma}, \cG) \times H^n(\bA_K^{\in\Sigma}, \cG\pr).
\]

\item \lab{top-gp}
For $n = 0$, and also for $n \ge 1$ if $\cG$ is commutative, $H^n(\bA_K^{\in\Sigma}, \cG)$ is a topological group.
\eenum
\eprop

\bpf 
We let $\cG_{\Sigma_0}$ and $\cG_{\Sigma_0}\pr$ be $\bA_{K, \Sigma_0}^{\in\Sigma}$-descents of $\cG$ and $\cG\pr$ as in \S\ref{ad-topo-Hn}.
\benum
\item
We enlarge $\Sigma_0$ to descend $\cG \ra \cG\pr$ to an $\bA_{K, \Sigma_0}^{\in\Sigma}$-homomorphism $\cG_{\Sigma_0} \ra \cG\pr_{\Sigma_0}$. Then it remains to apply \cite{Ces15d}*{3.3}, which supplies the continuity of 
\[
\qqq H^n(\cO_v, \cG_{\Sigma_0}) \ra H^n(\cO_v, \cG_{\Sigma_0}\pr) \qq \text{for $v\in \Sigma\setminus \Sigma_0$}
\]
and of 
\[
\qqq H^n(K_v, \cG) \ra H^n(K_v, \cG\pr) \qq \text{for $v\in \Sigma$.}
\]

\item
The bijections in question are continuous by \ref{basic-a}. Their openness reduces to that of similar bijections for the targets of \eqref{ad-topo-def}. The openness of \eqref{open-tr} reduces further to proving the openness of similar bijections for the sources of \eqref{open-tr}. For this, the openness of 
\[\ba
\qqq H^n(\cO_v, \cG_{\Sigma_0} \times \cG\pr_{\Sigma_0}) &\isomto H^n(\cO_v, \cG_{\Sigma_0}) \times H^n(\cO_v, \cG\pr_{\Sigma_0})\quad \text{ for $v\in \Sigma\setminus \Sigma_0$} \quad \quad \text{ and of} \\
\qqq H^n(K_v, \cG \times \cG\pr) &\isomto H^n(K_v, \cG) \times H^n(K_v, \cG\pr) \quad\quad\ \, \text{ for $v\in  \Sigma$}
\ea\]
provided by \cite{Ces15d}*{3.5~(c)~and~2.9~(b)} (with \Cref{Art-input} as in \S\ref{ad-topo-Hn}) suffices. 

\item
To check the continuity of the multiplication and the inverse of $H^n(\bA_K^{\in\Sigma}, \cG)$ we proceed similarly to the proof of \ref{basic-b} to reduce to the analogous continuity for $H^n(\cO_v, \cG_{\Sigma_0})$ with $v\in \Sigma\setminus \Sigma_0$ and for $H^n(K_v, \cG)$ with $v \in \Sigma$. To conclude we apply \cite{Ces15d}*{3.5 (c) and 3.6}.   \qedhere
\eenum
\epf

\bprop \lab{twist}
Let $\cG$ be a flat $\bA_K^{\in\Sigma}$-group algebraic space of finite presentation, $\cT$ a right $\cG$-torsor, and ${}_\cT \cG \ce \Aut_\cG(\cT)$ the resulting inner form of $\cG$. Twisting by $\cT$ induces a homeomorphism 
\[
H^1(\bA_K^{\in\Sigma}, {}_\cT \cG) \cong H^1(\bA_K^{\in\Sigma}, \cG).
\]
\eprop

\bpf
We enlarge $\Sigma_0$ to descend $\cT$ to a $\cG_{\Sigma_0}$-torsor $\cT_{\Sigma_0}$. For an $\bA_{K, \Sigma_0}^{\in\Sigma}$-algebra $R$, we consider the natural in $R$ twisting by $\cT_{\Sigma_0}$ bijection 
\[
H^1(R, {}_{\cT_{\Sigma_0}} \cG_{\Sigma_0}) \cong H^1(R, \cG_{\Sigma_0}),
\] 
which by \cite{Ces15d}*{3.4} is a homeomorphism for $R = \cO_v$ with $v\in \Sigma\setminus \Sigma_0$ and for $R = K_v$ with $v \in \Sigma$. Thus, twisting by $\cT_{\Sigma_0}$ induces a bijection
\[
\varinjlim_{\Sigma\pr_0} \p{ \prod_{v\in \Sigma\pr_0} H^1(K_v, {}_{\cT}\cG) \times \prod_{v\in \Sigma \setminus \Sigma\pr_0} H^1(\cO_v, {}_{\cT_{\Sigma_0}}\cG_{\Sigma_0}) } \cong \varinjlim_{\Sigma\pr_0}  \p{ \prod_{v\in \Sigma\pr_0} H^1(K_v, \cG) \times \prod_{v\in\Sigma \setminus \Sigma\pr_0} H^1(\cO_v, \cG_{\Sigma_0}) }
\]
that is a homeomorphism compatible with the maps \eqref{ad-pullb}. The claim follows.
\epf

\bprop\lab{disc}
For an $\bA_K^{\in\Sigma}$-group algebraic space $\cG$ of finite presentation, if
\begin{enumerate}[label=(\arabic*)]
\item \lab{disc-sm}
$n = 1$, and $\cG$ is smooth with connected fibers; or

\item \lab{disc-2}
$n \ge 2$, and $\cG$ is commutative and either smooth or isomorphic to $G_{\bA_K^{\in\Sigma}}$ for a $K$-group scheme~$G$,
\eenum
then $H^n(\bA_K^{\in\Sigma}, \cG)$ is discrete.
\eprop

\bpf
By \Cref{res-pr-top} \ref{res-pr-top-b} and \Cref{Hens-0} (with \Cref{limit}), 
\[
\tst H^n(\bA_K^{\in\Sigma}, \cG) \cong \bigoplus_{v\in \Sigma} H^n(K_v, \cG)
\] 
both algebraically and topologically. By \cite{Ces15d}*{3.5}, every singleton of $\bigoplus_{v\in \Sigma} H^n(K_v, \cG)$ is open.
\epf

We do not mention the $n \ge 2$ cases explicitly in \Cref{haus,lc,open-maps}---for such $n$ we have nothing to add to the discreteness that \Cref{disc} \ref{disc-2} offers.

\bprop \lab{haus}
For a flat $\bA_K^{\in \Sigma}$-group algebraic space $\cG$ of finite presentation, if 
\benuma
\item \lab{haus-0}
$n = 0$, and $\cG$ is separated; or

\item \lab{haus-2}
$n = 1$, and $\cG$ is either smooth with connected fibers or proper,
\eenum
then $H^n(\bA_K^{\in \Sigma}, \cG)$ is Hausdorff. 
\eprop

\bpf
By \Cref{H1-inj}, \eqref{ad-topo-def} is injective. By \Cref{res-pr-top} \ref{res-pr-top-a} and \cite{Ces15d}*{2.17~(3), 3.8~(1), and 3.8~(4)}, the target of \eqref{ad-topo-def} is Hausdorff.
\epf

\bprop \lab{lc}
For a flat $\bA_K^{\in \Sigma}$-group algebraic space $\cG$ of finite presentation, if 
\benuma
\item
$n = 0$, and $\cG$ is separated; or

\item
$n = 1$, and $\cG$ satisfies any of the following: is isomorphic to $G_{\bA_K^{\in\Sigma}}$ for a proper $K$-group scheme $G$, is smooth with connected fibers, or is smooth and proper,
\eenum
then $H^n(\bA_K^{\in \Sigma}, \cG)$ is locally compact.
\eprop

\bpf
By \Cref{res-pr-top} \ref{res-pr-top-b}, $H^n(\bA_K^{\in \Sigma}, \cG)$ is homeomorphic to $\prod\pr_{v\in\Sigma} H^n(K_v, \cG)$, which is locally compact by \cite{Ces15d}*{3.7 and 3.10} (with \Cref{Art-input}).
\epf

We conclude \S\ref{topo-ad} with a study of some topological aspects of long exact cohomology sequences. 

\bpp[Continuity of maps in cohomology sequences] \lab{seq-cont}
Let 
\[
\iota\colon \cH \hra \cG
\]
be a monomorphism of flat $\bA_K^{\in \Sigma}$-group algebraic spaces of finite presentation, and 
\[
\cQ \ce \cG/\iota(\cH)
\]
the resulting quotient left homogeneous space, which is a flat $\bA_K^{\in \Sigma}$-algebraic space of finite presentation by \cite{SP}*{\href{http://stacks.math.columbia.edu/tag/06PH}{06PH}, \href{http://stacks.math.columbia.edu/tag/06ET}{06ET}, \href{http://stacks.math.columbia.edu/tag/06EV}{06EV}, \href{http://stacks.math.columbia.edu/tag/040W}{040W}, and \href{http://stacks.math.columbia.edu/tag/03KS}{03KS}}. We consider the cohomology sequence
\be \lab{LES}  \tag{$\dagger$}
\dotsb \ra \cQ(\bA_K^{\in \Sigma}) \ra H^1(\bA_K^{\in \Sigma}, \cH) \ra H^1(\bA_K^{\in \Sigma}, \cG) \xra{x} H^1(\bA_K^{\in \Sigma}, \cQ) \xra{y} H^2(\bA_K^{\in \Sigma}, \cH) \xra{z} \dotsb,
\ee
in which $x$ (resp.,~$y$) is defined if $\iota(\cH)$ is normal (resp.,~central) in $\cG$, and $z$ and the subsequent maps are defined if $\cG$ is commutative. We use \cite{SP}*{\href{http://stacks.math.columbia.edu/tag/07SM}{07SM}} to enlarge a $\Sigma_0$ of \S\ref{ad-topo-Hn} and  descend $\iota$ to a monomorphism $\iota_{\Sigma_0} \colon \cH_{\Sigma_0} \hra \cG_{\Sigma_0}$ and $\cQ$ to $\cQ_{\Sigma_0} \ce \cG_{\Sigma_0}/\iota_{\Sigma_0}(\cH_{\Sigma_0})$; we enlarge $\Sigma_0$ further to ensure that $\iota_{\Sigma_0}(\cH_{\Sigma_0})$ is normal (resp.,~central) in $\cG_{\Sigma_0}$ if $\iota(\cH)$ is normal (resp.,~central) in $\cG$. By \cite{Ces15d}*{4.2}, all the maps of the counterparts of \eqref{LES} over $\cO_v$ for $v\in \Sigma \setminus \Sigma_0$ and over $K_v$ for $v\in \Sigma$ are continuous (whenever defined), and hence so are all the maps of \eqref{LES}.
\epp

\bprop \lab{open-maps}
Assume the setup of \S\ref{seq-cont}.
\benum
\item \lab{open-maps-a}
If $\cH$ is smooth with connected fibers, then the map 
\[
\cG(\bA_K^{\in \Sigma}) \ra \cQ(\bA_K^{\in \Sigma}) 
\]
is open.

\item \lab{open-maps-b}
If $\cG$ is smooth with connected fibers, and $\cH$ is either smooth or isomorphic to the base change of a $K$-group scheme, then the map 
\[
\cQ(\bA_K^{\in \Sigma}) \ra H^1(\bA_K^{\in \Sigma}, \cH) 
\]
is open.

\item \lab{open-maps-c}
If $\cQ$ is smooth with geometrically connected fibers, and each of $\cH$ and $\cG$ is either smooth or isomorphic to the base change of a $K$-group scheme, then the map 
\[
H^1(\bA_K^{\in \Sigma}, \cH) \ra H^1(\bA_K^{\in \Sigma}, \cG) 
\]
is open.

\item \lab{open-maps-d}
If $\cH$ is commutative, $\iota(\cH)$ is normal in $\cG$, and each of $\cG$ and $\cQ$ is either smooth or isomorphic to the base change of a $K$-group scheme, then the map 
\[
H^1(\bA_K^{\in \Sigma}, \cG) \xra{x} H^1(\bA_K^{\in \Sigma}, \cQ) 
\]
is open.
\eenum
\eprop

\bpf 
We apply \Cref{limit} to enlarge $\Sigma_0$ to ensure that $\bA_{K, \Sigma_0}^{\in \Sigma}$-descents inherit smoothness or fibral geometric connectedness. By \Cref{Bha-input} and \Cref{ad-sm-main,ad-main}, the $H^n(\bA_K^{\in \Sigma}, -)$ in \ref{open-maps-a}--\ref{open-maps-d} have bijective comparison maps \eqref{ad-topo-def}. Thus, by \cite{Ces15d}*{4.3} (i.e.,~by an analogue of \ref{open-maps-a}--\ref{open-maps-d} with $\bA_K^{\in \Sigma}$ replaced by $\cO_v$ with $v\in \Sigma\setminus \Sigma_0$ or by $K_v$ with $v\in \Sigma$) and by the openness of \eqref{open-tr}, it remains to argue the following claims:
\begin{itemize}
\item For \ref{open-maps-a}, the surjectivity of $\cG_{\Sigma_0}(\cO_v) \ra \cQ_{\Sigma_0}(\cO_v)$ for $v\in \Sigma\setminus \Sigma_0$;

\item For \ref{open-maps-b}, the surjectivity of $\cQ_{\Sigma_0}(\cO_v) \ra~H^1(\cO_v, \cH_{\Sigma_0})$ for $v\in \Sigma\setminus \Sigma_0$;

\item For \ref{open-maps-c}, the surjectivity of $H^1(\cO_v, \cH_{\Sigma_0}) \ra H^1(\cO_v, \cG_{\Sigma_0})$ for $v\in \Sigma\setminus \Sigma_0$;

\item For \ref{open-maps-d}, the surjectivity of $H^1(\cO_v, \cG_{\Sigma_0}) \ra H^1(\cO_v, \cQ_{\Sigma_0})$ for $v\in \Sigma\setminus \Sigma_0$.
\end{itemize}

We turn to proving these surjectivity claims.
\benum
\item[(a)--(b)]
The surjectivity follows by combining \Cref{Hens-0} with the cohomology sequence over $\cO_v$.

\addtocounter{enumi}{2}
\item
To lift a $(\cG_{\Sigma_0})_{\cO_v}$-torsor $\cT$ to an $(\cH_{\Sigma_0})_{\cO_v}$-torsor, we need to show that $\cT/(\cH_{\Sigma_0})_{\cO_v}$, which is fppf locally isomorphic to $(\cQ_{\Sigma_0})_{\cO_v}$, has an $\cO_v$-point. Since $\cT/(\cH_{\Sigma_0})_{\cO_v}$ is a left homogeneous space under the inner form $\Aut_{\cG_{\Sigma_0}}(\cT)$ of $\cG_{\Sigma_0}$, it has an $\bF_v$-point by \Cref{Lang-revamp} \ref{Lang-revamp-b}. To lift an $\bF_v$-point of the smooth $\cO_v$-algebraic space $\cT/(\cH_{\Sigma_0})_{\cO_v}$ to an $\cO_v$-point, we apply \cite{Ces15d}*{B.8}.

\item
Let $\bbB \cG_{\Sigma_0}$ and $\bbB \cQ_{\Sigma_0}$ be the classifying stacks of $\cG_{\Sigma_0}$ and $\cQ_{\Sigma_0}$. By \cite{Ces15d}*{A.4 (d) and B.8}, $\bbB \cG_{\Sigma_0} \ra \bbB \cQ_{\Sigma_0}$ is a smooth morphism of $\bA_{K, \Sigma_0}^{\in \Sigma}$-algebraic stacks whose essential surjectivity on $\cO_v$-points follows from that on $\bF_v$-points supplied by \Cref{Lang-revamp} \ref{ff-H1-surj}. \qedhere
\eenum
\epf


\section{Discreteness of the image of global cohomology inside adelic cohomology} \lab{disc-image}

The main goals of this section are \Cref{Sha} \ref{SF-d} and \Cref{ad-dual,im-closed}. These results go most of the way towards proving the sought Poitou--Tate sequence \eqref{PT-new}, which is derived in \S\ref{PT}. Topological arguments that involve the adelic topologies of \S\ref{topo-ad} play an important role throughout \S\S\ref{disc-image}--\ref{CPT}, so we begin the discussion by summarizing in \S\ref{basic-topo} the necessary topological facts. We proceed to use these facts without explicit notice.

\bpp[$H^n(\bA_K^{\in \Sigma}, G)$ for a commutative finite $G$] \lab{basic-topo}
Let $G$ be a commutative finite $K$-group scheme and let $\cG \ra U$ be its commutative finite flat model over a nonempty open $U \subset S$. For $n \ge 0$, we endow $H^n(\bA_K^{\in \Sigma}, G)$ with the topology of \S\ref{ad-topo-Hn}, so \Cref{res-pr-top} \ref{res-pr-top-b} provides a homeomorphism 
\[
\tst H^n(\bA_K^{\in \Sigma}, G) \cong \prod\pr_{v\in \Sigma} H^n(K_v, G)
\]
induced by the pullback maps, where the restricted product is with respect to the open subgroups 
\[
H^n(\cO_v, \cG) \subset H^n(K_v, G) \qq \text{with $v\in U \cap \Sigma$.}
\] 
By Propositions \ref{basic} \ref{top-gp}, \ref{disc}~\ref{disc-2}, \ref{haus}, and \ref{lc}, each $H^n(\bA_K^{\in \Sigma}, G)$ is a locally compact Hausdorff abelian topological group that is discrete if $n \ge 2$. If $n \ge 3$, then $H^n(\bA_K^{\in \Sigma}, G)$ is finite, and, in fact, then
\[
\tst H^n(\bA_K^{\in \Sigma}, G) \cong \prod_{v \text{ real}} H^n(K_v, G).
\]
\epp

\bpp[The $\wh{H}^0$ modification] \lab{H0-Tate}
Throughout \S\ref{disc-image}--\ref{CPT}, if $v$ is an archimedean place, then we implicitly make the Tate modification: $H^n(K_v, -)$ stands for $\wh{H}^n(K_v, -)$ (this has no effect for $n \ge 1$). We therefore write $\wh{H}^0(\bA_K^{\in \Sigma}, -)$ for $\prod\pr_{v \in \Sigma} H^0(K_v, -)$, where the Tate modification is in effect for archimedean $v$, and reserve $H^0(\bA_K^{\in \Sigma}, -)$ for honest global sections. In the setting of \S \ref{basic-topo}, 
\[
\tst H^0(\bA_K^{\in \Sigma}, G) \cong \prod_{v\in \Sigma} H^0(K_v, G),
\]
so there is a surjection 
\[
H^0(\bA_K^{\in \Sigma}, G) \ra \wh{H}^0(\bA_K^{\in \Sigma}, G).
\]

The principal advantage of $\wh{H}^0(\bA_K^{\in \Sigma}, -)$ over $H^0(\bA_K^{\in \Sigma}, -)$ is the perfectness of the pairings \eqref{ad-pair-H02}.
\epp

\bpp[The setup] \lab{setup}
Throughout \S\ref{disc-image}--\ref{CPT} we fix commutative Cartier dual finite flat $(S - \Sigma)$-group schemes $G$ and $H$ and adopt the following notation (together with its analogue for~$H$):
\benuma
\item
For $n\in \bZ_{\ge 0}$, we let 
\[
\loc^n(G)\colon H^n(S - \Sigma, G) \ra H^n(\bA_K^{\in \Sigma}, G)
\] 
denote the indicated pullback, where $H^n(\bA_K^{\in \Sigma}, G)$ stands for $\wh{H}^0(\bA_K^{\in \Sigma}, G)$ if $n = 0$; 
\item 
We set $\Sha^n(G) \ce \Ker(\loc^n(G))$, so that we also have
\be \lab{Sha-alt-def}
\textstyle \Sha^n(G) = \Ker\p{H^n(S - \Sigma, G) \ra \prod_{v\in \Sigma} H^n(K_v, G)}.
\ee
Likewise, for a commutative finite flat group scheme $\cG$ over a nonempty open $U \subset S$, we set
\[
\textstyle \Sha^n(\cG) \ce \Ker\p{H^n(U, \cG) \ra \prod_{v\not\in U} H^n(K_v, \cG)}.
\]
\eenum
\epp

\blem\lab{Sha}
We fix an $n \in \bZ_{\ge 0}$ and spread out $G$ and $H$ to Cartier dual commutative finite flat $U$-group schemes $\cG$ and $\cH$ for some open $U \subset S$ that contains $S - \Sigma$.
\benum
\item \lab{SF-a}
If $n \not\in \{0, 1, 2, 3\}$ (resp.,~if $U \neq S$ and $n \not\in \{1, 2\}$), then 
\[
\qq \Sha^n(\cG) = 0 \qq \text{and} \qq \Sha^n(\cH) = 0.
\] 
If $n \in \{0, 1, 2, 3\}$, then there is a perfect bilinear pairing of finite abelian groups
\be \lab{Sha-pair}
\Sha^n(\cG) \times \Sha^{3 - n}(\cH) \xra{\text{\cite{Ces14a}*{2.16}}} \bQ/\bZ.
\ee

\item \lab{SF-c}
If $U$ is small enough, then the pullback $H^n(U, \cG) \ra H^n(S - \Sigma, G)$ maps $\Sha^n(\cG)$ isomorphically onto $\Sha^n(G)$, and likewise with $H$ in place of $G$. Thus, $\Sha^n(G)$ and $\Sha^n(H)$ are~finite.

\item \lab{SF-d}
If $n \not\in \{0, 1, 2, 3\}$ (resp.,~if $(S - \Sigma) \neq S$ and $n \not\in \{1, 2\}$), then 
\[
\qq \Sha^n(G) = 0 \qq \text{and} \qq \Sha^n(H) = 0.
\] 
Moreover, for small enough $U$, \eqref{Sha-pair} induces perfect pairings
of finite abelian groups
\be \lab{Sha-pair-pair}
\qq\Sha^1(G) \times \Sha^2(H) \ra \bQ/\bZ \quad \quad \text{and} \quad \quad \Sha^2(G) \times \Sha^1(H) \ra \bQ/\bZ.
\ee
\eenum
\elem

Parts \ref{SF-a} and \ref{SF-c} of \Cref{Sha} serve as buildup to the key part \ref{SF-d}: part \ref{SF-a} proves \ref{SF-d} in the case when $\Sigma$ is finite, whereas \ref{SF-c} facilitates reduction to this special case.

\brem \lab{depend-on-U}
The pairings \eqref{Sha-pair-pair} a priori depend on the choice of a small enough $U$. Our reliance on the Artin--Verdier duality proved in \cite{Mil06}*{III.\S8} forces us to use flat cohomology with compact supports defined in \cite{Mil06}*{III.0.6 (b)} using completions, and not Henselizations as in \cite{Mil06}*{III.0.4}, which makes the functorialities needed for the independence of $U$ harder to check. 
\erem

\bpf[Proof of \Cref{Sha}]
\hfill
\benum
\item
The perfectness of \eqref{Sha-pair} was established in \cite{Ces14a}*{2.16} (and earlier in \cite{GA09}*{4.7} in a special case). Thus, the claimed vanishing for $n = 3$ follows from that for $n = 0$, which is evident from \eqref{Sha-alt-def}. The vanishing for $n \ge 4$ follows from \cite{Mil06}*{III.0.6 (a), III.3.2, and III.8.2}, according to which the map in \eqref{Sha-alt-def} is an isomorphism for $n \ge 4$ if $\Sigma$ is finite.

\item
We let $U\pr$ range over the open subsets of $U$ that contain $S - \Sigma$. 

\bcl \lab{unif-bdd}
For each $U\pr$, one has $\#\Sha^n(\cG_{U\pr}) \le \#\Sha^n(\cG)$.
\ecl

\bpf
Part \ref{SF-a} reduces to the $n = 1$ case, in which the Cartesian square
\be \ba \lab{H1-H1}
\xymatrix{
H^1(U, \cG) \ar@{^(->}[r] \ar[d] & H^1(U\pr, \cG) \ar[d] \\
\prod_{v\in U \setminus U\pr} H^1(\cO_v, \cG) \ar@{^(->}[r] & \prod_{v\in U \setminus U\pr} H^1(K_v, \cG)
}
\ea\ee
leads to the inclusion $\Sha^1(\cG_{U\pr}) \subset \Sha^1(\cG)$ (for a proof of Cartesianness, see \cite{Ces15b}*{4.3}).
\epf

By \Cref{lim}, 
\[
\tst H^n(S - \Sigma, G) = \varinjlim_{U\pr} H^n(U\pr, \cG_{U\pr}), 
\]
so a fixed $x \in \Sha^n(G)$ comes from $\Sha^n(\cG_{U\pr})$ for every small enough $U\pr$. Thus, \ref{SF-a} and \Cref{unif-bdd} prove the finiteness of $\Sha^n(G)$ and supply a $U\pr$ that works for every $x$. We shrink this $U\pr$ until $\#\Sha^n(\cG_{U\pr}) = \#\Sha^n(G)$ to arrange $\Sha^n(\cG_{U\pr}) \cong \Sha^n(G)$. We then rename the new $U\pr$ to $U$ and repeat the same steps with $G$ replaced by $H$ to arrive at a desired small enough $U$.

\item
The claim follows by combining \ref{SF-a} and \ref{SF-c}.
\qedhere
\eenum
\epf

\brems
\remi
\Cref{Sha} \ref{SF-c} implies in particular that fixed $x, y \in H^n(S - \Sigma, G)$ have the same pullback to every $H^n(K_v, G)$ with $v \in \Sigma$ as soon as such pullbacks agree for $v$ in some large finite subset of $\Sigma$. To see this, choose $U$ in \Cref{Sha} \ref{SF-c} so that $x$ and $y$ come from $H^n(U, \cG)$ and let the places outside $U$ comprise the subset in question.

\remi
See \cite{Mil06}*{I.4.10 (a)} and \cite{GA09}*{4.6--4.8} for special cases of \Cref{Sha}.
\erems

\bpp[The invariant map] 
We define it to be the composite 
\[
\textstyle \inv\colon H^2(\bA_K^{\in \Sigma}, \bG_m) \overset{\ref{ad-sm-main}}{\cong} \bigoplus_{v \in \Sigma} H^2(K_v, \bG_m) \xra{\Sigma_{v\in \Sigma} \inv_v} \bQ/\bZ,
\]
where the isomorphism is induced by the pullbacks and the $\inv_v$ are the local invariant maps.
\epp

\bpp[The duality pairings] \lab{pair}
Cup product and Cartier duality give the bilinear pairing
\be\lab{ad-pair-H1}
H^1(\bA_K^{\in \Sigma}, G) \times H^{1}(\bA_K^{\in \Sigma}, H) \ra H^2(\bA_K^{\in \Sigma}, \bG_m) \xra{\inv} \bQ/\bZ,
\ee
as well as the analogous pairings
\be\lab{ad-pair-H02-pre}
H^0(\bA_K^{\in \Sigma}, G) \times H^{2}(\bA_K^{\in \Sigma}, H) \ra \bQ/\bZ \quad\quad \text{and} \quad\quad H^2(\bA_K^{\in \Sigma}, G) \times H^{0}(\bA_K^{\in \Sigma}, H) \ra \bQ/\bZ.
\ee
The pairings \eqref{ad-pair-H1} and \eqref{ad-pair-H02-pre} identify with the sum of the pairings
\be \lab{local-dual}
H^n(K_v, G) \times H^{2 - n}(K_v, H) \ra H^2(K_v, \bG_m) \xra{\inv_v} \bQ/\bZ 	\quad\quad \quad \text{for $v \in \Sigma$},
\ee
in which the Tate modification may be taken 
to be in play for archimedean $v$. The pairings \eqref{ad-pair-H02-pre} therefore induce bilinear pairings
\be\lab{ad-pair-H02}
\wh{H}^0(\bA_K^{\in \Sigma}, G) \times H^{2}(\bA_K^{\in \Sigma}, H) \ra \bQ/\bZ \quad\quad \text{and} \quad\quad H^2(\bA_K^{\in \Sigma}, G) \times \wh{H}^{0}(\bA_K^{\in \Sigma}, H) \ra \bQ/\bZ.
\ee
\epp

\bprop \lab{ad-dual} \hfill
\benum
\item \lab{AD-a}
If $\bQ/\bZ$ is taken to be discrete, then the pairings \eqref{ad-pair-H1}--\eqref{ad-pair-H02} are continuous.

\item \lab{AD-b}
The pairings \eqref{ad-pair-H1}, \eqref{local-dual} with the Tate modification in play, and \eqref{ad-pair-H02} are perfect, i.e., they exhibit their terms as Pontryagin dual locally compact Hausdorff abelian topological~groups.

\item \lab{AD-c}
The images of the following pairs of maps are orthogonal complements under \eqref{ad-pair-H1} or~\eqref{ad-pair-H02}: 
\[\ba
H^0(S - \Sigma, G) \xra{\loc^0(G)} \wh{H}^0(\bA_K^{\in \Sigma}, G) \quad \quad \text{ and }\quad \quad H^{2}(S - \Sigma, H) \xra{\loc^2(H)} H^{2}(\bA_K^{\in \Sigma}, H),\\
H^1(S - \Sigma, G) \xra{\loc^1(G)} H^1(\bA_K^{\in \Sigma}, G) \quad \quad \text{ and }\quad \quad H^{1}(S - \Sigma, H) \xra{\loc^1(H)} H^{1}(\bA_K^{\in \Sigma}, H), \\
H^2(S - \Sigma, G) \xra{\loc^2(G)} H^2(\bA_K^{\in \Sigma}, G) \quad \quad \text{ and }\quad \quad H^{0}(S - \Sigma, H) \xra{\loc^0(H)} \wh{H}^{0}(\bA_K^{\in \Sigma}, H).
\ea\]
\eenum
\eprop

\emph{Proof.}
We spread out $G$ and $H$ to Cartier dual commutative finite flat $U$-group schemes $\cG$ and $\cH$ for an open $U \subset S$ that contains $S - \Sigma$.
\benum
\item
The invariant maps are continuous because their sources are discrete by \Cref{disc}~\ref{disc-2}. Moreover, the first map of \eqref{local-dual} is continuous by \cite{Ces14a}*{A.3}. It remains to note that if $v\in U \cap \Sigma$, then the pairing \eqref{local-dual} vanishes on $H^n(\cO_v, \cG) \times H^{2 - n}(\cO_v, \cH)$ (e.g.,~by~\Cref{Hens-0}). 

\item
Tate--Shatz local duality \cite{Mil06}*{I.2.3, I.2.13 (a), III.6.10} ensures perfectness of \eqref{local-dual}.\footnote{The agreement of the topologies used here and in loc.~cit.~is guaranteed by \cite{Ces15d}*{3.5 (c), 5.11, and 6.5}. }
Perfectness of \eqref{ad-pair-H1} and \eqref{ad-pair-H02} then follows from \cite{Tat50}*{Thm.~3.2.1}, because
\be \lab{o-comp}
\quad \quad H^n(\cO_v, \cG) \subset H^n(K_v, G) \quad \text{ and } \quad H^{2 - n}(\cO_v, \cH) \subset H^{2 - n}(K_v, H) \quad \quad \text{for $v\in U \cap \Sigma$}
\ee
are compact open orthogonal complements by \cite{Ces15d}*{3.10} and \cite{Mil06}*{III.1.4 and III.7.2}.

\item
In this proof, we abbreviate $\wh{H}^0(\bA_K^{\in \Sigma}, -)$ by $H^0(\bA_K^{\in\Sigma}, -)$. By \cite{Ces14c}*{5.3}, the images of
\be \lab{im-oc}
\textstyle 
\quad \quad H^n(U, \cG) \ra \bigoplus_{v\not\in U} H^n(K_v, G) \quad\quad \text{and} \quad\quad H^{2 - n}(U, \cH) \ra \bigoplus_{v\not\in U} H^{2 - n}(K_v, H)
\ee
are orthogonal complements. Therefore, by shrinking $U$ and combining the orthogonality of the subgroups in \eqref{o-comp} with the limit formalism supplied by \Cref{lim}, we deduce the orthogonality part of the claim. 

We use \Cref{Sha} \ref{SF-c} to shrink $U$ to arrange that the pullback induces an isomorphism $\Sha^n(\cG_{U\pr}) \cong \Sha^n(G)$ for every open $U\pr \subset U$ that contains $S - \Sigma$ and likewise for $H$. We let  
\[
\textstyle \qqq x = (x_v)_{v\in \Sigma}  \quad\quad \text{in} \quad\quad H^n(\bA_K^{\in \Sigma}, G) = \prod\pr_{v\in \Sigma} H^n(K_v, G)
\]
be orthogonal to $\im(\loc^{2 - n}(H))$ and we shrink $U$ to ensure that $x_v \in H^n(\cO_v, \cG)$ for $v \in U \cap \Sigma$. 

Due to the settled orthogonality part, the following claim together with its analogue for $H$ will finish the proof.

\bcl \lab{x-image}
The element $x$ lies in the image of 
\[
\qqq \loc^n(G)\colon H^n(S - \Sigma, G) \ra H^n(\bA_K^{\in\Sigma}, G).
\]
\ecl

\emph{Proof.}
The images of the maps in \eqref{im-oc} are orthogonal complements, so $(x_v)_{v \not\in U}$ is the pullback of some $y \in H^n(U, \cG)$. We replace $x$ by $x - y|_{H^n(\bA_K^{\in \Sigma}, G)}$ to reduce to the case when $x_v = 0$ for $v\not\in U$. Due to \Cref{Hens-0}, this settles the $n = 2$ case. 

In the $n = 1$ case, we let $U\pr \subset U$ be an open that contains $S - \Sigma$, so that \eqref{im-oc} also gives a $y\pr \in H^n(U\pr, \cG)$ that pulls back to $(x_v)_{v\not\in U\pr}$. Since the square \eqref{H1-H1} is Cartesian, $y\pr \in H^1(U, \cG)$, and even $y\pr \in \Sha^1(\cG)$. But then the choice of $U$ forces $y\pr \in \Sha^1(\cG_{U\pr})$, so $x_v = 0$ for $v\not\in U\pr$. We vary $U\pr$ to conclude that $x = 0$.

In the $n = 0$ case, we proceed similarly: if $\Sigma$ is finite, we choose $U = S - \Sigma$; else, we choose $U \neq S$, vary $U\pr$, and use the injectivity of 
\[
\qqq H^0(U\pr, \cG) \ra H^0(K_v, G) \qq  \text{for $v \nmid \infty$}
\]
to prove that $x = 0$.  \QEDD
\eenum

\bcor
If $n \ge 3$, then $\loc^n(G)$ is a surjection 
\[
\tst H^n(S - \Sigma, G) \surjects \prod_{\text{real } v} H^n(K_v, G).
\]
If either $n \ge 3$ with $(S - \Sigma) \neq S$ or $n \ge 4$, then $\loc^n(G)$ is an isomorphism 
\[
\tst H^n(S - \Sigma, G) \cong \prod_{\text{real } v} H^n(K_v, G).
\]
\ecor

\bpf
If $\#\Sigma < \infty$, then the fact that the images of \eqref{im-oc} are orthogonal complements (proved in \cite{Ces14c}*{5.3}) gives the surjectivity claim. Surjectivity in general then follows from a limit argument. To obtain the claim of the last sentence it remains to apply \Cref{Sha} \ref{SF-d}.
\epf

\bprop \lab{im-closed}
For $n \in \bZ_{\ge 0}$, the image 
\[
\im(\loc^n(G)) \subset H^n(\bA_K^{\in \Sigma}, G)
\] 
is closed, discrete, and cocompact (here and in the proof, $H^n$ for $n = 0$ abbreviates $\wh{H}^0$), and likewise with $H$ in place of $G$.
\eprop

\bpf
If $n \ge 3$, then $H^n(\bA_K^{\in \Sigma}, G)$ itself is finite and discrete, so we assume that $n \le 2$. 

\Cref{ad-dual} \ref{AD-a} and \ref{AD-c} imply closedness. \Cref{ad-dual} \ref{AD-c} and \cite{HR79}*{24.11} then identify 
\[
\tst H^n(\bA_K^{\in\Sigma}, G)/\im(\loc^n(G)) \qq \text{and}\qq  \im(\loc^{2 - n}(H))
\]
as Pontryagin duals, so it remains to prove discreteness.

If $n = 2$, then $H^n(\bA_K^{\in\Sigma}, G)$ itself is discrete. If $n = 0$, then the discreteness of $\im(\loc^n(G))$ follows from its finiteness. If $n = 1$, then we first spread out $G$ to a commutative finite flat $U$-group scheme $\cG$ for an open $U \subset S$ that contains $S - \Sigma$. Then we apply \cite{Ces15b}*{4.3} to deduce that 
\[
\tst \im(\loc^1(G)) \cap \p{\prod_{v\in U \cap \Sigma} H^1(\cO_v, \cG) \times \prod_{v\not\in U} H^1(K_v, G)}
\]
is the image of the pullback map
\[
\textstyle H^1(U, \cG) \ra \prod_{v\in U \cap \Sigma} H^1(\cO_v, \cG) \times \prod_{v\not\in U} H^1(K_v, G).
\]
The discreteness of $\im(\loc^1(G))$ now follows from \cite{Ces14a}*{2.9}, which ensures that the pullback 
\[
\tst H^1(U, \cG) \ra \bigoplus_{v\not \in U} H^1(K_v, G)
\]
has a discrete image and a finite kernel.
\epf

\begin{q} \lab{qq}
Is the image of the pullback map 
\[
H^1(S - \Sigma, G) \ra H^1(\bA_K^{\in \Sigma}, G)
\]
still closed and discrete if $G$ is only assumed to be a separated $(S - \Sigma)$-group scheme of finite type?
\end{q}

We do not know the answer to \Cref{qq} even in the case when $\Sigma$ consists of all the places.


\section{Poitou--Tate} \lab{PT}

With the results of \S\ref{disc-image} at hand, we are ready for the sought Poitou--Tate sequence \eqref{PT-seq}. We continue to assume the setup of \S\ref{setup} and adopt the convention introduced in \S\ref{H0-Tate}.

\bthm \lab{PT-thm}
For Cartier dual commutative finite flat $(S - \Sigma)$-group schemes $G$ and $H$, the~sequence
\be\ba \lab{PT-seq}
\xymatrix@C=36pt @R=12pt{
    &  H^0(S - \Sigma, G) \ar[r]^-{\loc^0(G)} & \wh{H}^0(\bA_K^{\in \Sigma}, G)  \ar[r]^-{\loc^2(H)^*} & H^2(S - \Sigma, H)^* 
                \ar@{->} `r/6pt[d] `/7pt[l] `^dl/12pt[lll] `^r/9pt[dl] [dll] \\
             & H^1(S - \Sigma, G) \ar[r]^-{\loc^1(G)} & H^1(\bA_K^{\in \Sigma}, G)  \ar[r]^-{\loc^1(H)^*} & H^1(S - \Sigma, H)^* 
                \ar@{->} `r/6pt[d] `/7pt[l] `^dl/12pt[lll] `^r/9pt[dl] [dll] \\
           &     H^2(S - \Sigma, G) \ar[r]^-{\loc^2(G)} & H^2(\bA_K^{\in \Sigma}, G)  \ar[r]^-{\loc^0(H)^*} & H^0(S - \Sigma, H)^* 
}
\ea\ee
is exact, where the $\loc^n(H)^*$ are defined using the perfectness of \eqref{ad-pair-H1} and \eqref{ad-pair-H02} supplied by \Cref{ad-dual} \ref{AD-b}, the curved arrows are defined using the perfectness of \eqref{Sha-pair-pair},\footnote{If $\Sigma$ is infinite, then the curved arrows inherit from \eqref{Sha-pair-pair} an a priori dependence on noncanonical choices, see \Cref{depend-on-U}.} and all the maps are continuous when $H^n(S - \Sigma, G)$ and $H^n(S - \Sigma, H)$ are endowed with the discrete topology. Moreover, if $\p{S - \Sigma} \neq S$, then $\loc^0(G)$ is injective and $\loc^0(H)^*$ is surjective.
\ethm

\bpf
In this proof, $H^n(\bA_K^{\in \Sigma}, -)$ for $n = 0$ abbreviates $\wh{H}^n(\bA_K^{\in \Sigma}, -)$.

For $n \in \{0, 1, 2\}$, \cite{HR79}*{24.11} with \Cref{ad-dual} identifies the exact sequences 
\be \lab{very-temp}\ba
0 \ra &\im(\loc^n(G)) \ra  H^n(\bA_K^{\in \Sigma}, G) \ra H^n(\bA_K^{\in\Sigma}, G)/\im(\loc^n(G)) \ra 0 \quad \quad \quad \quad  \text{and}\\
0 \ra &\p{H^{2 - n}(\bA_K^{\in\Sigma}, H)/\im(\loc^{2 - n}(H))}^* \ra  H^{2 - n}(\bA_K^{\in\Sigma}, H)^* \ra \im(\loc^{2 - n}(H))^* \ra 0.
\ea\ee
Moreover, \cite{HR79}*{24.11} with the discreteness aspect of \Cref{im-closed} gives the exact sequence
\be \lab{dual-ses}
0 \ra \im(\loc^{2- n}(H))^* \ra H^{2 - n}(S - \Sigma, H)^* \ra \Sha^{2 - n}(H)^* \ra 0.
\ee
The exactness of \eqref{PT-seq} and the continuity of its maps follow. The claim about the injectivity of $\loc^0(G)$ and the surjectivity of $\loc^0(H)^*$ follows from \Cref{Sha} \ref{SF-d} and \eqref{dual-ses}.
\epf

\brems
\remi \lab{strict}
The maps in \eqref{PT-seq} are strict, i.e.,~the quotient and subspace topologies on their images agree: for $\loc^n(G)$, this follows from \Cref{im-closed}; for $\loc^n(H)^*$, this follows from \eqref{very-temp} and \eqref{dual-ses}; for the curved arrows, one also uses the finiteness of $\Sha^{2 - n}(H)^*$ supplied by \Cref{Sha} \ref{SF-d}. 

\remi
If $(S - \Sigma) = S$ and either $K$ has no real places or $\#G$ is odd, then \eqref{PT-seq} recovers dualities between 
\[
H^1(S, G) \qq \text{and} \qq H^2(S, H)
\]
and between 
\[
H^2(S, G) \qq \text{and}\qq H^1(S, H).
\]
\erems

\beg
We present a concrete cohomological computation that illustrates the expectation that a suitable version of the original Poitou--Tate sequence \eqref{PT-old} may hold without the assumption on the ``ramification set'' $\Sigma$ (\Cref{PT-thm} fulfills this expectation). For this, we will compute the order of $H^2(S - \Sigma, \mu_p)$ in terms of other invariants of $K$ and $\Sigma$, will observe that the computation does not need assumptions on the residue characteristics of the finite places in $\Sigma$, and will then show that the same order may also be computed from the Poitou--Tate sequence with $G = \bZ/p\bZ$ and $H = \mu_p$. 

To simplify the computation we will let $K$ be a number field, $p$ be an odd prime, and $\Sigma$ be a finite set of places of $K$ that contains at least one nonarchimedean place and all the archimedean places. 

From the cohomology of the exact sequence $0 \ra \mu_p \ra \bG_m \xra{p} \bG_m \ra 0$ over $S - \Sigma$ we get
\[
0 \ra \Pic(S - \Sigma)/p\Pic(S - \Sigma) \ra H^2(S - \Sigma, \mu_p) \ra \Br(S - \Sigma)[p] \ra 0.
\]
From the residue sequence \cite{Gro68c}*{2.1--2.2} for the Brauer group and from the global class field theory sequence of \Cref{CFT}, we deduce the short exact sequence
\[
\tst 0 \ra \Br(S - \Sigma)[p] \ra \bigoplus_{v\in \Sigma} \Br(K_v)[p] \ra \bZ/p\bZ \ra 0.
\]
Since $p$ is odd, this sequence gives the formula
\[
\#\Br(S - \Sigma)[p] = p^{\#(\Sigma \cap \{v \nmid \infty\}) - 1},
\]
which leads to the sought expression for the order of $H^2(S - \Sigma, \mu_p)$:
\be \lab{H2-card}
\#H^2(S - \Sigma, \mu_p) = p^{\#(\Sigma \cap \{v \nmid \infty\}) - 1} \cdot \#\Pic(S - \Sigma)[p].
\ee
As promised, the computation needed no assumptions on residue characteristics of finite places in $\Sigma$.

Let us now recover the same formula \eqref{H2-card} from the Poitou--Tate sequence \eqref{PT-seq}. With the choice $G = \bZ/p\bZ$ and $H = \mu_p$, the Poitou--Tate sequence gives the exact sequence
\be\ba \lab{PT-take}
\xymatrix@C=16pt @R=12pt{
& & 0  \ra  \bZ/p\bZ \ra \bigoplus_{v \in \Sigma, \, v \nmid \infty} \bZ/p\bZ    \ra H^2(S - \Sigma, \mu_p)^*
                \ar@{->} `r/7pt[d] `/7pt[l] `^dl/12pt[l] `^r/9pt[d] [d] & \\
   &        & \Ker\p{\Hom(\pi_1^\et(S - \Sigma), \bZ/p\bZ) \ra \bigoplus_{v \in \Sigma} \Hom(\Gal(\ov{K}_v/K_v), \bZ/p\bZ)}   \ar[r] & 0.
}
\ea\ee
The displayed kernel is dual to the Galois group of the maximal abelian unramified extension of $K$ in which all the places of $\Sigma$ split and whose Galois group is killed by $p$. Therefore, by global class field theory, this kernel has the same cardinality as $\Pic(S - \Sigma)/p \Pic(S - \Sigma)$, and from \eqref{PT-take} we recover the formula \eqref{H2-card}. However, if we were to use the version of \eqref{PT-take} that results from the original Poitou--Tate sequence \eqref{PT-old}, then in this second computation of $\#H^2(S - \Sigma, \mu_p)$ we would have had to assume that $\Sigma$ contains all the places of $K$ above $p$.
\eeg


\section{Cassels--Poitou--Tate} \lab{CPT}

The Cassels--Poitou--Tate sequence alluded to in the section title is presented in \Cref{CPT-thm}. It relates certain Selmer subgroups of $H^1(S - \Sigma, G)$ to those of $H^1(S - \Sigma, H)$. For its proof, which is based on the Poitou--Tate sequence \eqref{PT-seq}, topological considerations facilitated by the results of \S\S\ref{topo-ad}--\ref{disc-image} are key even in the number field case because $H^1(\bA_K^{\in \Sigma}, G)$ need not be discrete when $\Sigma$ is infinite. In the case of finite $\Sigma$, \Cref{CPT-thm} was proved in \cite{Ces14a}*{4.2}. Although the proof below goes along the same lines, it does not seem possible to directly reduce to this special case.

Throughout \S\ref{CPT} we continue to assume the setup of \S\ref{setup}.

\bpp[Local conditions that are orthogonal complements] \lab{sel}
Suppose that we have subgroups
\[
\Sel(G_{\bA_K^{\in \Sigma}}) \subset H^1(\bA_K^{\in\Sigma}, G) \quad \quad \text{and} \quad \quad \Sel(H_{\bA_K^{\in \Sigma}}) \subset H^1(\bA_K^{\in\Sigma}, H)
\]
that are orthogonal complements under the perfect pairing \eqref{ad-pair-H1} (cf.~\Cref{ad-dual} \ref{AD-b}). Define the resulting Selmer groups $\Sel(G)$ and $\Sel(H)$ by the exactness of
\[
\ba
&0 \ra \Sel(G) \ra H^1(S - \Sigma, G) \ra H^1(\bA_K^{\in\Sigma}, G)/\Sel(G_{\bA_K^{\in \Sigma}}) \quad \quad \quad \quad \text{and} \\
&0 \ra \Sel(H) \ra H^1(S - \Sigma, H) \ra H^1(\bA_K^{\in\Sigma}, H)/\Sel(H_{\bA_K^{\in \Sigma}}).
\ea
\]
\epp

\bthm \lab{CPT-thm}
In the setup of \S\ref{sel}, suppose that each of $\Sel(G_{\bA_K^{\in \Sigma}})$ and $\Sel(H_{\bA_K^{\in \Sigma}})$ is either open or compact (when endowed with the subspace topology). Then the sequence
\be\ba\lab{CPT-seq}
\xymatrix@C=36pt @R=12pt{
  0  \ar[r] &   \Sel(G) \ar[r] & H^1(S - \Sigma, G)  \ar[r] & \f{H^1(\bA_K^{\in \Sigma}, G)}{\Sel(G_{\bA_K^{\in\Sigma}})} \ar[r]^{y(G)} & \Sel(H)^*
                \ar@{->} `r/9pt[d] `/9pt[l] `^dl/12pt[lll] `^r/9pt[dl] [dll]^-{x(G)} \\
             & & H^2(S - \Sigma, G) \ar[r]^-{\loc^2(G)} & H^2(\bA_K^{\in\Sigma}, G)  \ar[r]^-{\loc^0(H)^*} & H^0(S - \Sigma, H)^* 
}
\ea\ee
is exact, where $y(G)$ is the dual of 
\[
\textstyle \Sel(H) \xra{\loc^1(H)|_{\Sel(H)}} \Sel(H_{\bA_K^{\in \Sigma}})
\]
and $x(G)$ is the map factoring the second curved arrow of \eqref{PT-seq}.\footnote{In particular, if $\Sigma$ is infinite, then $x(G)$ inherits an a priori dependence on noncanonical choices.} If $\Sel(G)$, $\Sel(H)$, $H^n(S - \Sigma, G)$, and $H^0(S - \Sigma, H)$ are taken to be discrete, then the maps in \eqref{CPT-seq} are continuous. If $(S - \Sigma) \neq S$, then $\loc^0(H)^*$ is surjective.
\ethm

\bpf
We begin by justifying the definitions of $y(G)$ and $x(G)$. That of $y(G)$ rests on the isomorphism 
\[
\tst \f{H^1(\bA_K^{\in\Sigma}, G)}{\Sel(G_{\bA_K^{\in \Sigma}})} \cong \Sel(H_{\bA_K^{\in\Sigma}})^*
\]
obtained from \eqref{ad-pair-H1} by combining \Cref{ad-dual} \ref{AD-b} with \cite{HR79}*{24.11}. Loc.~cit.~together with the combination of \Cref{im-closed} and \cite{BouTG}*{III.28, Cor.~3} (the latter two are not needed if $\Sel(H_{\bA_K^{\in\Sigma}})$ is open) also gives the exact sequences in the commutative diagram
\be\ba \lab{big-com}
\xymatrix@C=16pt @R=16pt{
0 \ar[r] & \p{\f{H^1(\bA_K^{\in\Sigma}, H)}{\Sel(H_{\bA_K^{\in \Sigma}})}}^* \ar@{->>}[d] \ar[r] & H^1(\bA_K^{\in\Sigma}, H)^* \ar[d] \ar[r] &\Sel(H_{\bA_K^{\in \Sigma}})^* \ar[d]^-{y(G)} \ar[r] & 0, \\
0 \ar[r] &\im\p{H^1(S - \Sigma, H) \ra \f{H^1(\bA_K^{\in\Sigma}, H)}{\Sel(H_{\bA_K^{\in \Sigma}})}}^* \ar[r] & H^1(S - \Sigma, H)^* \ar[r] &\Sel(H)^* \ar[r] &0.
}
\ea\ee
The indicated surjectivity, and hence also the existence of the claimed $x(G)$, follows by also using 

\bcl \lab{cl-aux}
The following subgroups are orthogonal complements (and hence are closed):
\[
\im(\loc^1(G)|_{\Sel(G)}) \subset \Sel(G_{\bA_K^{\in\Sigma}}) \quad \text{and}  \quad \im\p{H^1(S - \Sigma, H) \ra \f{H^1(\bA_K^{\in\Sigma}, H)}{\Sel(H_{\bA_K^{\in \Sigma}})}} \subset \f{H^1(\bA_K^{\in\Sigma}, H)}{\Sel(H_{\bA_K^{\in \Sigma}})}.
\]
\ecl

\bpf
Since $\im(\loc^1(H)) + \Sel(H_{\bA_K^{\in \Sigma}})$ is closed in $H^1(\bA_K^{\in\Sigma}, H)$ (by \cite{BouTG}*{III.28, Cor.~1} if $\Sel(H_{\bA_K^{\in \Sigma}})$ is compact), \Cref{ad-dual} \ref{AD-c} and \cite{HR79}*{24.10} provide orthogonal complements
\[
\im(\loc^1(G)|_{\Sel(G)}) \subset H^1(\bA_K^{\in\Sigma}, G) \quad \quad \text{and} \quad \quad \im(\loc^1(H)) + \Sel(H_{\bA_K^{\in \Sigma}}) \subset H^1(\bA_K^{\in\Sigma}, H),
\]
which give rise to the claimed orthogonal complements.
\epf

Continuity of the maps in \eqref{CPT-seq} follows from their definitions. Surjectivity of $\loc^0(H)^*$ when $(S - \Sigma) \neq S$ follows from \Cref{PT-thm}. We turn to the remaining exactness claim. 
\begin{itemize}
\item
Exactness at $\Sel(G)$ and at $H^1(S - \Sigma, G)$ amounts to the definition of $\Sel(G)$. 

\item
Exactness at $\f{H^1(\bA_K^{\in\Sigma}, G)}{\Sel(G_{\bA_K^{\in \Sigma}})}$ follows from \cite{HR79}*{24.11}, the analogue of \Cref{cl-aux} for $G$, and the discreteness aspect of \Cref{im-closed} (see the proof of \Cref{PT-thm} for an analogous argument). 

\item
Exactness at $\Sel(H)^*$ follows from the definition of $x(G)$ and the commutativity of \eqref{big-com}. 

\item
Exactness at $H^2(S - \Sigma, G)$ and at $H^2(\bA_K^{\in\Sigma}, G)$ follows from the exactness of \eqref{PT-seq}. \qedhere
\end{itemize} 
\epf

\brems
\remi
Similarly to Remark \ref{strict}, the maps in \eqref{CPT-seq} are strict: for all the maps except $y(G)$ and $x(G)$, this follows from Remark \ref{strict} and \cite{BouTG}*{III.28, Cor.~3}; for $x(G)$, one also uses the bottom row of \eqref{big-com}; for $y(G)$, this follows from the discreteness of the image of $\loc^1(H)|_{\Sel(H)}$ supplied by \Cref{im-closed}. 

\remi
With the choice 
\[
\qq \Sel(G_{\bA_K^{\in\Sigma}}) = 0 \qq  \text{and}  \qq \Sel(H_{\bA_K^{\in\Sigma}}) = H^1(\bA_K^{\in \Sigma}, H),
\]
the sequence \eqref{CPT-seq} specializes to a segment of the Poitou--Tate sequence \eqref{PT-seq}. The choice 
\[
\Sel(G_{\bA_K^{\in\Sigma}}) = H^1(\bA_K^{\in \Sigma}, G) \qq \text{and} \qq \Sel(H_{\bA_K^{\in\Sigma}}) = 0
\] 
recovers the duality between $\Sha^2(G)$ and $\Sha^1(H)$ proved in \Cref{Sha} \ref{SF-d}.

\remi
See \cite{Cre12}*{Prop.~4.6} for another proof of the exactness of \eqref{CPT-seq} at the $\f{H^1(\bA_K^{\in \Sigma}, G)}{\Sel(G_{\bA_K^{\in\Sigma}})}$ term in the case when $K$ is a number field, $\Sigma$ is the set of all places, and $\Sel(H_{\bA_K^{\in\Sigma}}) =\prod_{v\in T} U_v \times \prod_{v\not\in T} \{0\}$ with $T$ a finite subset of $\Sigma$ and $U_v$ a subgroup of $H^1(K_v, H)$ for $v\in T$.
\erems


\appendix

\section{The sequence \eqref{PT-new} recovers the sequence \eqref{PT-old}} \lab{app}

We prove that \eqref{PT-old} is indeed a special case of \eqref{PT-new}, as promised in \S\ref{PT-wo}. We begin with a well-known \Cref{ff-ner-mod}. For its statement we use \cite{BLR90}*{1.2/1} as our definition of a N\'{e}ron model.

\bprop\lab{ff-ner-mod} 
Let $D$ be a connected (and hence nonempty) Dedekind scheme, $k$ its function~field.
\benum
\item \lab{ff-0}
A finite \'{e}tale $D$-scheme $X$ is a N\'{e}ron model of its generic fiber.
\eenum
Let $F$ be a finite \'{e}tale $k$-group scheme and $k^s$ a separable closure of $k$.
\benum \addtocounter{enumi}{1}
\item\lab{ff-ner-mod-a} 
The N\'{e}ron model $\cF \ra D$ of $F$ exists and is separated quasi-finite \'{e}tale. It is finite if and only if $F(k^s)$ is unramified at all the nongeneric $d \in D$ (i.e., if and only if $G_{k_{d}^\sh}$ is constant for all such $d$, where $k_{d}^\sh$ is the fraction field of the strict Henselization of the local ring $\cO_{D, d}$).

\item \lab{ff-ner-mod-d}
The functor $\cG \mapsto \cG_k$ is an equivalence between the category of \'{e}tale group N\'{e}ron models and that of finite \'{e}tale $k$-group schemes. 

\item \lab{ff-ner-mod-e}
Commutative finite \'{e}tale $D$-group schemes $\cG$ form a full abelian subcategory of the category of abelian sheaves on $D_\et$. This category is equivalent by the functor $\cG \mapsto \cG(k^s)$ to the category of finite discrete $\Gal(k^s/k)$-modules that are unramified at all the nongeneric points of $D$.
\eenum
\eprop

\bpf \hfill
\benum
\item
In checking the N\'{e}ron property of $X$, \'{e}tale descent reduces to the case of a constant $X$. 

\item
For the existence, \ref{ff-0}, spreading-out, and \cite{BLR90}*{1.4/1 and 6.5/3} reduce to the case of a strictly local $D$, in which, by \cite{BLR90}*{7.1/1}, $\cF \ra D$ is obtained from $F$ by extending the constant subgroup $\underline{F(k)}_k \subset F$ to a constant subgroup over $D$. The other claims are immediate from this construction. 

\item
Essential surjectivity follows from \ref{ff-ner-mod-a}. Full faithfulness follows from the N\'{e}ron property. 

\item 
For the abelian subcategory claim, it suffices to note that the subcategory is closed under direct sums, kernels, and quotients. The rest follows from \ref{ff-0}, \ref{ff-ner-mod-a}, and \ref{ff-ner-mod-d}.
\qedhere
\eenum
\epf

\Cref{ff-ner-mod} \ref{ff-ner-mod-e} allows us to identify the module $M$ of \eqref{PT-old} with a commutative finite \'{e}tale $(S - \Sigma)$-group scheme $\cM$. Since $\#\cM$ is a unit on $(S - \Sigma)$, \cite{TO70}*{p.~17, Lemma 5} applies and proves that the Cartier dual $\cM^D$ of $\cM$ is $(S - \Sigma)$-\'{e}tale. Thus, $\cM^D(k^s)$ may be identified with $M^D$. To conclude that the sequence \eqref{PT-new} applied to $G = \cM$ and $H = \cM^D$ recovers \eqref{PT-old}, it remains to note: 
\benuma
\item \lab{1}
One has 
\[
\qq H^n(S - \Sigma, \cM) \cong H^n(\Gal(K_\Sigma/K), M),
\]
and similarly for $\cM^D$;

\item \lab{2}
One has 
\[
\qq\textstyle  H^n(\bA_K^{\in \Sigma}, \cM) \cong \prod_{v\in \Sigma}\pr H^n(K_v, M),
\]
where the restricted product is taken with respect to the unramified cohomology subgroups and the Tate modifications are in play for $n = 0$.
\eenum

\bpf[Proof of \ref{1}]
If $\Sigma$ is finite, then a claimed isomorphism is provided by \cite{Mil06}*{II.2.9}. To reduce to this case, we spread out $\cM$ to a finite \'{e}tale $(S - \Sigma_0)$-group scheme $\wt{\cM}$ for some finite $\Sigma_0 \subset \Sigma$ for which $\#M$ is a unit on $S - \Sigma_0$ and combine the isomorphism 
\[
\textstyle H^n(\Gal(K_{\Sigma}/K), M) \cong \varinjlim_{\Sigma_0 \subset \Sigma_0\pr \subset \Sigma, \  \#\Sigma_0\pr < \infty} H^n(\Gal(K_{\Sigma_0\pr}/K), M)
\]
provided by \cite{Ser02}*{I.\S2, Prop.~8} with the analogous isomorphism provided by \Cref{lim}.
\epf

\bpf[Proof of \ref{2}]
Let $\Sigma_0$ and $\wt{\cM}$ be as in the proof of \ref{1}, so \Cref{res-pr-comp} identifies $H^n(\bA_K^{\in \Sigma}, \cM)$ with the restricted product $\prod_{v\in \Sigma}\pr H^n(K_v, \cM)$ with respect to the subgroups 
\[
H^n(\cO_v, \wt{\cM}) \subset H^n(K_v, \wt{\cM}) \qqq \text{for $v\in \Sigma \setminus \Sigma_0$.}
\]
Moreover, \cite{Gro68c}*{11.7~$1^{\circ})$} identifies $H^n(K_v, \cM)$ with Galois cohomology $H^n(K_v, M)$. It remains to argue that this identification restricts to an isomorphism 
\[
H^n(\cO_v, \wt{\cM}) \cong H^n_\nr(K_v, M) \qqq \text{for $v \in \Sigma \setminus \Sigma_0$.}
\]
For $n = 0$, this follows from the definitions and from the equality $\wt{\cM}(\cO_v) = \wt{\cM}(K_v)$. For $n \ge 2$, both sides vanish due to \Cref{Hens-0} and \cite{Ser02}*{II.\S5.5, Prop.~18 (c)}. For $n = 1$, one combines \Cref{ff-ner-mod} \ref{ff-0} with \cite{Ces15b}*{A.3 and A.4}. 
\epf

\end{document}